\RequirePackage{etex}
\documentclass[a4paper]{article}
\usepackage{amsmath,amssymb,amsthm}
\usepackage[top=3cm,bottom=3cm,left=2.5cm,right=2.5cm]{geometry}
\usepackage{microtype}
\usepackage{tikz}
\usetikzlibrary{arrows,cd}
\usepackage[hypertexnames=false]{hyperref}
\hypersetup{
	colorlinks=true,
	citecolor=red,
	linkcolor=blue,
	urlcolor=orange
	}
\usepackage{pdflscape}
\usepackage{subfigure}
\usepackage{cleveref}
\usepackage{autonum}
\usepackage{enumitem}
\setlist[enumerate,1]{label=(\arabic*),ref=(\arabic*)}
\setlist[enumerate,2]{label=(\alph*),ref=\theenumi. (\alph*)}

\numberwithin{equation}{section}
\theoremstyle{definition}

\newtheorem{theorem}{Theorem}[section]
\crefname{theorem}{theorem}{theorems}

\crefname{maintheorem}{Main Theorem}{Main Theorems}

\newtheorem{lemma}[theorem]{Lemma}
\crefname{lemma}{lemma}{lemmas}

\newtheorem{proposition}[theorem]{Proposition}
\crefname{proposition}{proposition}{propositions}

\newtheorem{corollary}[theorem]{Corollary}
\crefname{corollary}{corollary}{corollaries}

\crefname{definition}{definition}{definitions}

\crefname{definitionproposition}{Definition-Proposition}{Definition-Proposirions}

\newtheorem{example}[theorem]{Example}
\crefname{example}{example}{examples}

\newtheorem{remark}[theorem]{Remark}
\crefname{remark}{remark}{remarks}

\newtheorem{question}[theorem]{Question}
\crefname{question}{question}{questions}
\newcommand{\setmid}{\; \middle|\;}
\newcommand{\lmod}{\operatorname{\mathrm{\hspace{-2pt}-mod}}}

\newcommand{\Ext}{\operatorname{Ext}\nolimits}

\newcommand{\Hom}{\operatorname{Hom}\nolimits}

\newcommand{\induc}{{\operatorname{Ind}\nolimits}}
\newcommand{\restr}{{\operatorname{Res}\nolimits}}
\newcommand{\add}{\operatorname{\mathrm{add}}}
\newcommand{\stautilt}{\operatorname{\mathrm{s\tau-tilt}}}

\newcommand{\inertiagp}{I}

\newcommand{\Addresses}{{
  \bigskip
  \footnotesize

  Ryotaro~KOSHIO\, %\orcidlink{https://orcid.org/0000-0002-6753-562X}
  \par\nopagebreak
  \textsc{Department of Mathematics, Tokyo University of Science}
  \par\nopagebreak
	1-3, Kagurazaka, Shinjuku-ku, Tokyo, 162-8601, Japan
  \par\nopagebreak
  E-mail: \href{mailto:1120702@ed.tus.ac.jp}{1120702@ed.tus.ac.jp}

  \medskip

  Yuta~KOZAKAI\, %\orcidlink{https://orcid.org/0000-0001-6832-0148}
  \par\nopagebreak
  \textsc{Department of Mathematics, Tokyo University of Science}
	\par\nopagebreak
	1-3, Kagurazaka, Shinjuku-ku, Tokyo, 162-8601, Japan
	\par\nopagebreak
  E-mail: \href{mailto:kozakai@rs.tus.ac.jp}{kozakai@rs.tus.ac.jp}
}}
\title{
	Normal subgroups and support $\tau$-tilting modules
  \footnote{\emph{Mathematics Subject Classification} (2020). 20C20, 16G10.}
  \footnote{\emph{Keywords.} Support $\tau$-tilting modules, blocks of finite groups, induction functors, restriction functors}
	}
\author{Ryotaro~KOSHIO \and Yuta~KOZAKAI}
\date{\today}
\begin{document}
\maketitle
\begin{abstract}
  Let $\tilde{G}$ be a finite group, $G$ a normal subgroup of $\tilde{G}$ and $k$ an algebraically closed field of characteristic $p>0$. The first main result in this paper is to show that support $\tau$-tilting $k\tilde{G}$-modules satisfying some properties are support $\tau$-tilting modules as $kG$-modules too. As the second main result, we  give equivalent conditions for support $\tau$-tilting $k\tilde{G}$-modules to satisfy the above properties, and show that the set of the support $\tau$-tilting $k\tilde{G}$-modules with the properties is isomorphic to the set of $\tilde{G}$-invariant support $\tau$-tilting $kG$-modules as partially ordered sets. As an application, we show that the set of $\tilde{G}$-invariant support $\tau$-tilting $kG$-modules is isomorphic to the set of support $\tau$-tilting $k\tilde{G}$-modules in the case that the index $G$ in $\tilde{G}$ is a $p$-power. As a further application, we give a feature of vertices of indecomposable $\tau$-rigid  $k\tilde{G}$-modules. Finally, we give the block versions of the above results.
\end{abstract}

\section{Introduction}\label{introduction}
Since 2014 when \(\tau\)-tilting theory was introduced by T.~Adachi, O.~Iyama, and I.~Reiten \cite{MR3187626}, the theory continues to develop rapidly.
The main theme of the theory is the study of support \(\tau\)-tilting modules, and many researchers have given the work on these.
In fact, the support \(\tau\)-tilting modules over finite dimensional algebras are under the one-to-one correspondences with the various representation-theoretically important objects including two-term silting complexes \cite{MR3187626}, functorially finite torsion classes \cite{MR3187626}, left-finite semibricks \cite{MR4139031}, two-term simple-minded collections \cite{MR4139031,MR3178243}, and so on.
% Therefore, the theory is useful for the classification of two-term tilting complexes over group algebras or block algebras of finite groups.
In particular, the theory is expected to be helpful in solving Brou\'{e}'s abelian defect group conjecture because the theory is useful for the classification of two-term tilting complexes over group algebras or block algebras of finite groups.
Even though the studies on the \(\tau\)-tilting theory related to the modular representation theory of finite groups are very important for these reasons, there are few such studies.
Therefore, the authors have given the studies combining the above two theories
\cite{https://doi.org/10.48550/arxiv.2208.14680,KK2,MR4243358}.
All of them show that the induction functors from \(kG\)-modules to \(k\tilde{G}\)-modules give the poset homomorphisms from the support \(\tau\)-tilting modules over \(kG\) to those over \(k\tilde{G}\) under appropriate assumptions,
where \(G\) is a normal subgroup of a finite group \(\tilde{G}\) and \(k\) an algebraically closed field of characteristic \(p>0\).
Naturally, we are interested in the following \namecref{question group algebra}.
\begin{question}\label{question group algebra}
  When the restriction functor from \(k\tilde{G}\)-modules to \(kG\)-modules give the maps from the support \(\tau\)-tilting modules over \(k\tilde{G}\) to those over \(kG\)?
\end{question}

In regarding this \namecref{question group algebra}, in \cite{https://doi.org/10.48550/arxiv.2209.02992}, S.~Breaz, A.~Marcus, and G.~C.~Modoi gave a positive answer in case that the quotient group \(\tilde{G}/G\) is a \(p\)-prime group (i.e. the prime number \(p\) does not divide the order of the factor group \(\tilde{G}/G\)).
Therefore, we consider the case that \(\tilde{G}/G\) is not necessarily a \(p\)-prime group, and get the following positive answer for the \namecref{question group algebra}.

\begin{theorem}[\Cref{main-thm-1} and \Cref{cor order 2022-12-26 18:27:13}]\label{main1}
  Let \(\tilde{G}\) be a finite group, \(G\) a normal subgroup of \(\tilde{G}\),
  \(\tilde{M}\) a relatively \(G\)-projective support \(\tau\)-tilting \(k\tilde{G}\)-module,
  and \((\tilde{M}, \tilde{P})\) a corresponding support \(\tau\)-tilting pair.
  If it holds that \(\induc_{G}^{\tilde{G}}\restr_{G}^{\tilde{G}}\tilde{M}\in \add \tilde{M}\),
  then \(\restr_{G}^{\tilde{G}}\tilde{M}\) is a support \(\tau\)-tilting \(kG\)-module,
  and \((\restr_{G}^{\tilde{G}}\tilde{M}, \restr_{G}^{\tilde{G}}\tilde{P})\) a corresponding support \(\tau\)-tilting pair.
  Moreover, for relatively \(G\)-projective support \(\tau\)-tilting \(k\tilde{G}\)-modules \(\tilde{M}_1\) and \(\tilde{M}_2\) with the property that \(\induc_{G}^{\tilde{G}}\restr_{G}^{\tilde{G}}\tilde{M}_i\in \add \tilde{M}_i\) for \(i=1,2\), if \(\tilde{M}_1\geq\tilde{M}_2\) in \(\stautilt k\tilde{G}\), then \(\restr_{G}^{\tilde{G}}\tilde{M}_1 \geq \restr_{G}^{\tilde{G}}\tilde{M}_2\) in \(\stautilt kG\).
\end{theorem}

Let \(M\) be a \(\tilde{G}\)-invariant support \(\tau\)-tilting \(kG\)-module.
The first author showed that \(\induc_{G}^{\tilde{G}} M\) is a
support \(\tau\)-tilting \(k\tilde{G}\)-module
\cite[Theorem 3.2]{https://doi.org/10.48550/arxiv.2208.14680}.
We are interested in what is the image of
the set of \(\tilde{G}\)-invariant support \(\tau\)-tilting \(kG\)-modules
under the map induced by the induction functor \(\induc_{G}^{\tilde{G}}\).
Therefore, we give equivalent conditions to the assumption of \Cref{main1}, and finally we clarify what the image of the map induced by \(\induc_{G}^{\tilde{G}}\) is in the following \namecref{main2}.

\begin{theorem}[\Cref{main-thm-2} and \Cref{cor1}]\label{main2}
  Let \(\tilde{M}\) be a support \(\tau\)-tilting \(k\tilde{G}\)-module.
  Then the following conditions are equivalent:
  \begin{enumerate}
    \item \(\tilde{M}=_{\add} \induc_{G}^{\tilde{G}} M\) for some \(\tilde{G}\)-invariant
          support \(\tau\)-tilting \(kG\)-module \(M\).
    \item \(\induc_{G}^{\tilde{G}}\restr_{G}^{\tilde{G}} \tilde{M} \in \add \tilde{M}\) and \(\tilde{M}\)
          is relatively \(G\)-projective.
    \item \(S\otimes_k \tilde{M} \in \add \tilde{M}\)
          for each simple \(k(\tilde{G}/G)\)-module \(S\).
  \end{enumerate}
  Moreover, denoting by
  \((\stautilt kG)^{\tilde{G}}\) the subset of
  \(\stautilt kG\) consisting of \(\tilde{G}\)-invariant support \(\tau\)-tilting \(kG\)-modules and by \((\stautilt k\tilde{G})^{\star}\) the subset of \(\stautilt k\tilde{G}\) consisting of support \(\tau\)-tilting \(k\tilde{G}\)-modules satisfying the above equivalent conditions, the induction functor \(\induc_{G}^{\tilde{G}}\) induces a poset isomorphism
  \begin{equation}
    \begin{tikzcd}[row sep=0pt]
      (\stautilt kG)^{\tilde{G}} \ar[r,"\sim"]& (\stautilt k\tilde{G})^{\star}\\
      M\ar[r,mapsto]&\induc_{G}^{\tilde{G}} M.
    \end{tikzcd}
  \end{equation}
\end{theorem}

The studies on the vertices of indecomposable modules over group algebras have been done by many researchers for a long time, for example, see \cite{MR3318541, MR2422312, MR2845567, MR131454, MR554380, MR4092934}.
On the other hand, \(\tau\)-rigid modules over finite dimensional algebras are important classes
and there are many studies on the modules. They have nice properties and are essential objects for the representation theory, for example, see \cite{MR3187626,MR617088,MR3910476,MR3856858,MR3428959}.
Therefore, one of our interests is to give a feature of the vertices of indecomposable \(\tau\)-rigid modules over group algebras.
As a further application of \Cref{main2}, we give a feature of vertices of indecomposable \(\tau\)-rigid modules.

\begin{theorem}[See \Cref{tau-rigid}]
  Let \(\tilde{G}\) be a finite group.
  Then any indecomposable \(\tau\)-rigid \(k\tilde{G}\)-module has a vertex contained in a Sylow \(p\)-subgroup of \(\tilde{G}\) properly if and only if \(\tilde{G}\) has a normal subgroup of \(p\)-power index in \(\tilde{G}\).
\end{theorem}

As a natural question, we wonder if we get the block version of our theorems for group algebras.
In particular, we are interested in how we give the block versions of \Cref{main1,main2}.
As the results, we get the block versions of the \namecrefs{main1}.
Let \(G\) be a normal subgroup of a finite group \(\tilde{G}\) and \(B\) a block of \(kG\).
We denote by \(\inertiagp_{\tilde{G}}(B):=\left\{ \tilde{g}\in \tilde{G} \setmid \tilde{g}B\tilde{g}^{-1}=B \right\}\) the inertial group of \(B\) in \(\tilde{G}\).
\begin{theorem}[{See \Cref{main-thm-1 block general}}]
  \label{main-thm-1 block general intro}
  Let \(G\) be a normal subgroup of a finite group \(\tilde{G}\), \(B\) a block of \(kG\), \(\tilde{B}\) a block of \(k\tilde{G}\) covering \(B\), \(\beta\) the block of \(k\inertiagp_{\tilde{G}}(B)\) satisfying
  \begin{equation}
    \sum_{x\in [\tilde{G}/\inertiagp_{\tilde{G}}(B)]}x1_\beta x^{-1}=1_{\tilde{B}}
  \end{equation}
  and \(\tilde{M}\) a support \(\tau\)-tilting \(\tilde{B}\)-module.
  If it holds that \(\beta\induc_{G}^{\inertiagp_{\tilde{G}}(B)}\restr_{G}^{\inertiagp_{\tilde{G}}(B)}\beta\restr_{\inertiagp_{\tilde{G}}(B)}^{\tilde{G}}\tilde{M} \in \add \beta\restr_{\inertiagp_{\tilde{G}}(B)}^{\tilde{G}}\tilde{M}\)
  and \(\beta\restr_{\inertiagp_{\tilde{G}}(B)}^{\tilde{G}}\tilde{M}\) is relatively \(G\)-projective,
  then we have that \(\restr_{G}^{\inertiagp_{\tilde{G}}(B)}\beta\restr_{\inertiagp_{\tilde{G}}(B)}^{\tilde{G}}\tilde{M}\) is a support \(\tau\)-tilting \(B\)-module.
  Moreover, if \((\tilde{M}, \tilde{P})\) is a support \(\tau\)-tilting pair for \(\tilde{B}\)
  corresponding to \(\tilde{M}\), then the pair
  \begin{equation}
    (\restr_{G}^{\inertiagp_{\tilde{G}}(B)}\beta\restr_{\inertiagp_{\tilde{G}}(B)}^{\tilde{G}}\tilde{M}, \restr_{G}^{\inertiagp_{\tilde{G}}(B)}\beta\restr_{\inertiagp_{\tilde{G}}(B)}^{\tilde{G}}\tilde{P})
  \end{equation}
  is a support \(\tau\)-tilting pair for \(B\)
  corresponding to \(\restr_{G}^{\inertiagp_{\tilde{G}}(B)}\beta\restr_{\inertiagp_{\tilde{G}}(B)}^{\tilde{G}}\tilde{M}\).
\end{theorem}

\begin{theorem}[{See \Cref{main-thm-2 block general} and \Cref{cor1 block general}}]\label{main-thm-2 block general intro}
  Let \(G\) be a normal subgroup of a finite group \(\tilde{G}\), \(B\) a block of \(kG\), \(\tilde{B}\) a block of \(k\tilde{G}\) covering \(B\), \(\beta\) the block of \(k\inertiagp_{\tilde{G}}(B)\) satisfying
  \begin{equation}
    \sum_{x\in [\tilde{G}/\inertiagp_{\tilde{G}}(B)]}x1_\beta x^{-1}=1_{\tilde{B}}
  \end{equation}
  and \(\tilde{M}\) a support \(\tau\)-tilting \(\tilde{B}\)-module.
  Then the following conditions are equivalent:
  \begin{enumerate}
    \item \(\tilde{M}=_{\add} \tilde{B}\induc_{G}^{\tilde{G}} M\) for some \(\inertiagp_{\tilde{G}}(B)\)-invariant support \(\tau\)-tilting \(B\)-module \(M\).\label{main-thm-2 block general item1 intro}
    \item \(\beta\induc_{G}^{\inertiagp_{\tilde{G}}(B)}\restr_{G}^{\inertiagp_{\tilde{G}}(B)}\beta\restr_{\inertiagp_{\tilde{G}}(B)}^{\tilde{G}}\tilde{M} \in \add \beta\restr_{\inertiagp_{\tilde{G}}(B)}^{\tilde{G}}\tilde{M}\) and \(\beta\restr_{\inertiagp_{\tilde{G}}(B)}^{\tilde{G}}\tilde{M}\) is relatively \(G\)-projective.\label{main-thm-2 block general item2 intro}
    \item \(\beta(S\otimes_k \beta\restr_{\inertiagp_{\tilde{G}}(B)}^{\tilde{G}}\tilde{M}) \in \add \beta\restr_{\inertiagp_{\tilde{G}}(B)}^{\tilde{G}}\tilde{M}\) for each simple \(k(\inertiagp_{\tilde{G}}(B)/G)\)-module \(S\).\label{main-thm-2 block general item3 intro}
  \end{enumerate}
  Moreover, denoting by
  \((\stautilt B)^{\inertiagp_{\tilde{G}}(B)}\) the subset of
  \(\stautilt B\) consisting of \(\inertiagp_{\tilde{G}}(B)\)-invariant support \(\tau\)-tilting \(B\)-modules and by \((\stautilt \tilde{B})^{\star \star \star}\) the subset of \(\stautilt \tilde{B}\) consisting of support \(\tau\)-tilting \(\tilde{B}\)-modules satisfying the above equivalent conditions, the functor \(\tilde{B}\induc_{G}^{\tilde{G}}\) induces a poset isomorphism
  \begin{equation}
    \begin{tikzcd}[row sep=0pt]
      (\stautilt B)^{\inertiagp_{\tilde{G}}(B)} \ar[r,"\sim"]& (\stautilt \tilde{B})^{\star \star \star}\\
      M\ar[r,mapsto]&\tilde{B}\induc_{G}^{\tilde{G}} M.
    \end{tikzcd}
  \end{equation}
\end{theorem}

Our particular interest is the case that the index of \(G\) in \(\tilde{G}\) is a \(p\)-power.
In fact, under some assumptions, it is expected that
tilting complexes over the block \(B\) of \(kG\) give
those over the unique block \(\tilde{B}\) of \(k\tilde{G}\) covering \(B\) (for example, see \cite{MR2592757,https://doi.org/10.48550/arxiv.2207.12668, MR1889341}).
In this regard, the authors showed the following result in \cite{MR4243358}.

\begin{theorem}[{\cite[Theorem 1.2]{MR4243358}}]\label{KK1-main}
  Let \(G\) be a normal subgroup of a finite group \(\tilde{G}\)
  of \(p\)-power index in \(\tilde{G}\),
  \(B\) a block of \(kG\),
  and \(\tilde{B}\) the unique block of \(k\tilde{G}\) covering \(B\).
  Assume that the following two conditions are satisfied:
  \begin{enumerate}
    \item Any indecomposable \(B\)-module is \(I_{\tilde{G}}(B)\)-invariant.
    \item The set of isomorphism classes of basic support \(\tau\)-tilting \(B\)-modules is a finite set.
  \end{enumerate}
  Then the induction functor \(\induc_{G}^{\tilde{G}}\) induces an isomorphism from \(\stautilt B\)
  to \(\stautilt \tilde{B}\) of partially ordered sets.
\end{theorem}

This \namecref{KK1-main} can be applied to the case that the block \(B\) has a cyclic defect group,
but the two conditions limit the scope of its use.
For example, the theorem cannot be applied to the case
that \(p=2\), \(G\) is the alternating group \(A_4\) of degree \(4\)
and that \(\tilde{G}\) is the symmetric group \(S_4\) of degree \(4\), because the nontrivial simple \(kA_4\)-modules are not \(S_4\)-invariant.
Indeed, \(\stautilt kA_4\) is not isomorphic to \(\stautilt kS_4\)
because the number of isomorphism classes of simple \(kA_4\)-modules is three and that of \(kS_4\) is two.
However, we wonder if the induction functor might give some kinds of good relation between the special subsets of the two, and finally, as an application of \Cref{main-thm-2 block general intro},
we could get the following \namecref{p-extension block general intro}  which can be applied to the case of \(kA_4\) and \(kS_4\).
The following \namecref{p-extension block general intro} is a significant generalization of \Cref{KK1-main} and enables us to explain the phenomenon occurred in \cite[Example 3.9]{https://doi.org/10.48550/arxiv.2208.14680} (see \Cref{A4S4Example 2023-01-04 20:09:08}).

\begin{theorem}[{See \Cref{p-extension block general}}]\label{p-extension block general intro}
  Let \(G\) be a normal subgroup of a finite group \(\tilde{G}\), \(B\) a block of \(kG\) and \(\tilde{B}\) a block of \(k\tilde{G}\) covering \(B\).
  If the quotient group \(\inertiagp_{\tilde{G}}(B)/G\) is a
  \(p\)-group, then the functor \(\induc_{G}^{\tilde{G}}\) induces an isomorphism as partially ordered sets between \((\stautilt B)^{\inertiagp_{\tilde{G}}(B)}\) and \(\stautilt \tilde{B}\),
  where \((\stautilt B)^{\inertiagp_{\tilde{G}}(B)}\) is
  the subset of \(\stautilt B\) consisting of \(\inertiagp_{\tilde{G}}(B)\)-invariant support \(\tau\)-tilting \(B\)-modules.
\end{theorem}
% -----------------------------------
% Notation
% -----------------------------------
Throughout this paper, we fix the following notation:

Let \(k\) be an algebraically closed field of characteristic \(p>0\).
An algebra means a \(k\)-algebra.
For a finite dimensional algebra \(\Lambda\), a \(\Lambda\)-module means a finite dimensional left \(\Lambda\)-module.
For a \(\Lambda\)-module \(M\), we denote the Auslander-Reiten translate of \(M\) by \(\tau M\).
In case that \(\Lambda\) is a symmetric algebra,
\(\tau M\) is isomorphic to \(\Omega^2 M\).
We denote the category of all direct summands of finite direct sums of copies of \(M\) by \(\add M\).
For \(\Lambda\)-modules \(M\) and \(N\),
% we write \(M\leq_{\add} N\) if \(\add M \subset \add N\) and write \(M=_{\add}N\) if \(\add M=\add N\).
we write \(M=_{\add}N\) if \(\add M=\add N\).
% These relations are a partial order and an equivalence relation, respectively.
This relation is an equivalence relation.
We denote by \(\stautilt \Lambda\) the set of equivalence classes of support \(\tau\)-tilting \(\Lambda\)-modules under the equivalence relation \(=_{\add}\).
% Moreover, we write \(M\leq_{\add} N\) if \(\add M \subset \add N\)

Let \(G\) be a finite group and \(H\) a subgroup of \(G\).
We denote the restriction functor from \(kG\)-modules
to \(kH\)-modules by \({\rm Res}_H^G\)
and the induction functor \(kG\otimes_{kH}-\) from \(kH\)-modules
to \(kG\)-modules by \({\rm Ind}_H^G\).
We denote the trivial \(kG\)-module by \(k_G\).

Let \(\tilde{G}\) be a finite group,
\(G\) a normal subgroup of \(\tilde{G}\).
We denote a set of coset representatives of \(G\) in \(\tilde{G}\) by \([\tilde{G}/G]\).
For a \(kG\)-module\(M\) and \(\tilde{g} \in \tilde{G}\), we define a \(kG\)-module \(\tilde{g}M\) consisting of symbols \(\tilde{g}m\) as a set,
where \(m\in M\), and its \(kG\)-module structure is given by \(\tilde{g}m+\tilde{g}m':=\tilde{g}(m+m'), g(\tilde{g}m):=\tilde{g}(\tilde{g}^{-1}g\tilde{g}m)\) and \(\lambda(\tilde{g}m)=\tilde{g}(\lambda m)\) for \(m, m'\in M, g\in G\) and \(\lambda\in k\).
For a \(kG\)-module \(M\), we say that \(M\) is
\(\tilde{G}\)-invariant if \(M\) is isomorphic to \(\tilde{g}M\) for any \(\tilde{g}\in\tilde{G}\).

\section{Preliminaries}\label{sec preliminaries}
In this \namecref{sec preliminaries}, we give elementary facts on the modular representation theory which are helpful to prove our results.
\begin{proposition}[{See \cite[Lemma 8.5, Lemma 8.6]{MR860771}}]\label{Theorem: Frobenius and projective}
  Let \(G\) be a finite group, \(K\) a subgroup of \(G\), \(H\) a subgroup of \(K\).
  For any \(kG\)-module \(U\) and \(kH\)-module \(V\),  the following hold:
  \begin{enumerate}
    \item \(\restr^G_H U \cong \restr^K_H\restr^G_K U\).
    \item \(\induc^G_H V \cong \induc^G_K\induc^K_H V\).\label{transitive induc 2022-06-09 02:03:23}
    \item \(\induc_H^G(V\otimes_k \restr^G_H U)\cong (\induc_H^G V)\otimes_k U\).
    \item \(\Hom_{kG}(U, \induc_H^G V)\cong \Hom_{kH}(\restr_H^G U, V)\).\label{frob adjoint 2022-12-26 19:17:16}
    \item  \(\Hom_{kG}(\induc_H^G V, U) \cong\Hom_{kH}(V, \restr_H^G U)\).\label{frob adjoint 2022-06-02 15:44:10}
    \item The functors \(\restr_H^G\) and \(\induc_H^G\) send free modules (projective modules) to free modules (projective modules, respectively).\label{proj to proj 2022-06-02 12:44:53}
  \end{enumerate}
\end{proposition}
In the modular representation theory of finite groups, Mackey's decomposition formula is well-known and important.
We recall Mackey's decomposition formula for normal subgroups.
\begin{proposition}[{See \cite[Lemma 8.7]{MR860771}}]\label{mackey}
  Let \(G\) be a normal subgroup of a finite group \(\tilde{G}\)
  and \(M\) a \(kG\)-module.
  Then the following isomorphism as \(kG\)-modules holds:
  \begin{equation}
    \restr_{G}^{\tilde{G}}\induc_{G}^{\tilde{G}} M \cong \bigoplus_{x\in [\tilde{G}/G]} xM.
  \end{equation}
\end{proposition}

The following is known as Eckmann-Shapiro Lemma.

\begin{lemma}[{See \cite[Proposition 2.20.7]{MR3821517}}]\label{Eckmann-Shapiro}
  Let \(H\) be a finite group of a finite group \(G\),
  \(M\) a \(kH\)-module and \(N\) a \(kG\)-module.
  Then for all \(n\in \mathbb{N}\),
  there exists an isomorphism of \(k\)-vector spaces:
  \begin{equation}
    \Ext^n_{kH}(M, {\rm Res}^G_H N) \cong \Ext^n_{kG}({\rm Ind}^G_H M, N)
  \end{equation}
\end{lemma}

The following \namecref{commute} is a refinement of \cite[Lemma 3.1]{https://doi.org/10.48550/arxiv.2208.14680}
which requires the \(\tilde{G}\)-invariance for the \(kG\)-module.

\begin{lemma}\label{commute}
  Let \(G\) be a normal subgroup of \(\tilde{G}\)
  and \(M\) a \(kG\)-module.
  Then the following hold:
  \begin{enumerate}
    \item \(\induc_{G}^{\tilde{G}} (\Omega M) \cong \Omega(\induc_{G}^{\tilde{G}} M)\).\label{commute item 0}
    \item \(\induc_{G}^{\tilde{G}} (\tau M) \cong \tau(\induc_{G}^{\tilde{G}} M)\).
  \end{enumerate}
\end{lemma}

\begin{proof}
  We enough to show that the statement \ref{commute item 0} holds since \(\tau \cong \Omega^2\) for symmetric algebras.
  % We consider the following exact sequence of \(kG\)-module:
  % \begin{equation}
  %   \begin{tikzcd}
  %     0\ar[r]& \Omega M \ar[r] &\ar[r] P(M)& \ar[r] M \ar[r]& 0.
  %   \end{tikzcd}
  % \end{equation}
  % By applying the induction functor \(\induc_{G}^{\tilde{G}}\) to this exact sequence,
  % we have the following exact sequence of \(k\tilde{G}\)-module:
  % \begin{equation}
  %   \begin{tikzcd}
  %     0\ar[r] & \induc_{G}^{\tilde{G}}(\Omega M) \ar[r] & \induc_{G}^{\tilde{G}} P(M) \ar[r] & \induc_{G}^{\tilde{G}} M \ar[r] & 0.
  %   \end{tikzcd}
  % \end{equation}
  % This exact sequence has the following one as a direct summand clearly:
  % \begin{equation}
  %   \begin{tikzcd}
  %     0\ar[r]& \Omega (\induc_{G}^{\tilde{G}} M) \ar[r]& P(\induc_{G}^{\tilde{G}} M) \ar[r]& \induc_{G}^{\tilde{G}} M \ar[r]& 0.
  %   \end{tikzcd}
  % \end{equation}
  % Hence, there exists a projective \(kG\)-module \(Q\)
  There exists a projective \(kG\)-module \(Q\)
  such that \(\induc_{G}^{\tilde{G}}(\Omega M) \cong \Omega (\induc_{G}^{\tilde{G}} M)\oplus Q \)
  and that \(\induc_{G}^{\tilde{G}} P( M)\cong P(\induc_{G}^{\tilde{G}} M)\oplus Q\).
  Hence, we have that
  \begin{equation}
    \restr_{G}^{\tilde{G}}\induc_{G}^{\tilde{G}}(\Omega M) \cong \restr_{G}^{\tilde{G}}\Omega (\induc_{G}^{\tilde{G}} M)\oplus \restr_{G}^{\tilde{G}} Q,
  \end{equation}
  and the left-hand side is isomorphic to \(\bigoplus_{x\in [\tilde{G}/G]}  x\Omega M\)
  by \Cref{mackey}.
  However, each \(x\Omega M\) has no projective summands
  and the restricted module \(\restr_{G}^{\tilde{G}} Q\) is a projective \(kG\)-module by
  \Cref{Theorem: Frobenius and projective} \ref{proj to proj 2022-06-02 12:44:53},
  which implies that \(Q=0\).
  Therefore, we conclude that \(\induc_{G}^{\tilde{G}}(\Omega M) \cong \Omega (\induc_{G}^{\tilde{G}} M)\).
\end{proof}

\begin{lemma}\label{restriction-invariance}
  Let \(G\) be a normal subgroup of a finite group \(\tilde{G}\)
  and \(\tilde{M}\) be a \(k\tilde{G}\)-module.
  Then \(\restr_{G}^{\tilde{G}}\tilde{M}\) is a \(\tilde{G}\)-invariant \(kG\)-module.
\end{lemma}
\begin{proof}

  Take \(\tilde{g} \in \tilde{G}\) arbitrarily.
  We consider the map
  \begin{equation}
    \begin{tikzcd}[row sep=1pt]
      f\colon \restr_{G}^{\tilde{G}}\tilde{M} \ar[r]& \tilde{g}\restr_{G}^{\tilde{G}}\tilde{M}\\
      m\ar[r,mapsto]& \tilde{g}m.
    \end{tikzcd}
  \end{equation}
  Clearly, this map is linear and bijective.
  We only show that the map is \(kG\)-homomorphism,
  but for any \(g\in G\) and \(m\in \restr_{G}^{\tilde{G}}\tilde{M}\),
  it holds that
  \begin{equation}
    f(gm)=\tilde{g}gm = \tilde{g}g\tilde{g}^{-1}\tilde{g}m
    =g\cdot\tilde{g}m=g\cdot f(m).
  \end{equation}
\end{proof}

\section{Main Theorems}\label{sec main block}

In this section, we give theorems stated in Section \ref{introduction}
and their proofs.
Throughout this section, \(\tilde{G}\) means a finite group and \(G\) a normal subgroup of \(\tilde{G}\).

First, we start with a consideration on restricted modules of rigid modules and \(\tau\)-rigid modules.
Let \(\Lambda\) be a finite dimensional algebra.
We recall that a \(\Lambda\)-module \(M\) is rigid (resp. \(\tau\)-rigid)
if \(\Ext_{\Lambda}^{1}(M,M)=0\) (resp. \(\Hom_\Lambda(M, \tau M)=0\)).
We remark that \(\tau\)-rigid modules are rigid modules by Auslander-Reiten duality
\(\overline{\Hom}_\Lambda(X, Y)\cong D\Ext_{\Lambda}^1(Y, \tau X)\).

\begin{lemma}\label{res-rigid}
  Let \(\tilde{M}\) be a \(k\tilde{G}\)-module
  with the property that
  \(\induc_{G}^{\tilde{G}}\restr_{G}^{\tilde{G}}\tilde{M} \in \add \tilde{M}\).
  Then the following hold:
  \begin{enumerate}
    \item If \(\tilde{M}\) is a rigid \(k\tilde{G}\)-module, then the restricted module \(\restr_{G}^{\tilde{G}} \tilde{M}\) is a rigid \(kG\)-module.\label{res-rigid item1}
    \item If \(\tilde{M}\) is a \(\tau\)-rigid \(k\tilde{G}\)-module, then the restricted module \(\restr_{G}^{\tilde{G}} \tilde{M}\) is a \(\tau\)-rigid \(kG\)-module.\label{res-rigid item2}
  \end{enumerate}
\end{lemma}
\begin{proof}
  \ref{res-rigid item1}
  Let \(\tilde{M}\) be a rigid \(k\tilde{G}\)-module.
  Then, by \Cref{Eckmann-Shapiro},
  we have that
  \begin{equation}
    \Ext_{kG}^{1}(\restr_{G}^{\tilde{G}}\tilde{M}, \restr_{G}^{\tilde{G}}\tilde{M})
    \cong
    \Ext_{k\tilde{G}}^{1}(\tilde{M}, \induc_{G}^{\tilde{G}}\restr_{G}^{\tilde{G}}\tilde{M}).
  \end{equation}
  By the assumption that \(\induc_{G}^{\tilde{G}}\restr_{G}^{\tilde{G}}\tilde{M} \in \add \tilde{M}\)
  and the rigidity of \(\tilde{M}\),
  we have that the right-hand side is \(0\).
  Hence, \(\restr_{G}^{\tilde{G}}\tilde{M}\) is a rigid \(kG\)-module.

  \noindent
  \ref{res-rigid item2}
  Let \(\tilde{M}\) be a \(\tau\)-rigid \(k\tilde{G}\)-module.
  Then we have that
  \begin{equation}
    \Hom_{kG}(\restr_{G}^{\tilde{G}}\tilde{M}, \tau\restr_{G}^{\tilde{G}}\tilde{M})
    \cong
    \Hom_{k\tilde{G}}(\tilde{M}, \induc_{G}^{\tilde{G}}\tau\restr_{G}^{\tilde{G}}\tilde{M})
    \cong
    \Hom_{k\tilde{G}}(\tilde{M}, \tau\induc_{G}^{\tilde{G}}\restr_{G}^{\tilde{G}}\tilde{M}),
  \end{equation}
  where the last isomorphism comes from \Cref{commute}.
  By the assumption that \(\induc_{G}^{\tilde{G}}\restr_{G}^{\tilde{G}}\tilde{M} \in \add \tilde{M}\)
  and the \(\tau\)-rigidity of \(\tilde{M}\),
  we have that
  \(\Hom_{k\tilde{G}}(\tilde{M}, \tau\induc_{G}^{\tilde{G}}\restr_{G}^{\tilde{G}}\tilde{M})=0\),
  which implies that \(\restr_{G}^{\tilde{G}}\tilde{M}\) is a \(\tau\)-rigid \(kG\)-module.
\end{proof}

For a finite group \(H\) and a subgroup \(K\) of \(H\),
we recall that a \(kH\)-module \(M\) is relatively \(K\)-projective
if \(M\) is a direct summand of \(\induc_K^H\restr_K^H M\).

\begin{lemma}\label{tau-restriction}
  Let \(\tilde{M}\) be a relatively \(G\)-projective \(k\tilde{G}\)-module.
  Then \(\restr_{G}^{\tilde{G}}(\Omega\tilde{M})\cong \Omega(\restr_{G}^{\tilde{G}}\tilde{M})\).
  In particular, it holds that \(\tau(\restr_{G}^{\tilde{G}}\tilde{M})\cong\restr_{G}^{\tilde{G}}(\tau\tilde{M})\).
\end{lemma}
\begin{proof}
  There exists a projective \(kG\)-module \(P\)
  such that \(\restr_{G}^{\tilde{G}}(\Omega\tilde{M})\cong \Omega(\restr_{G}^{\tilde{G}}\tilde{M})\oplus P\).
  Hence, we enough to show that \(P=0\).
  It is clear in the case that \(\tilde{M}\) is a projective \(k\tilde{G}\)-module.

  We may assume that \(\tilde{M}\) has no projective summands.
  Since \(\tilde{M}\) is relatively \(G\)-projective, \(\Omega\tilde{M}\) is relatively \(G\)-projective too
  (for example see \cite[Proposition 20.7]{MR860771}).
  Hence, \(\Omega\tilde{M}\) is a direct summand of
  \(\induc_{G}^{\tilde{G}}\restr_{G}^{\tilde{G}}\Omega\tilde{M}\).
  On the other hand, by the isomorphism \(\restr_{G}^{\tilde{G}}\Omega\tilde{M}\cong \Omega(\restr_{G}^{\tilde{G}}\tilde{M})\oplus P\),
  we have that \(\induc_{G}^{\tilde{G}}\restr_{G}^{\tilde{G}}\Omega\tilde{M}\cong \induc_{G}^{\tilde{G}}(\Omega(\restr_{G}^{\tilde{G}}\tilde{M}))\oplus \induc_{G}^{\tilde{G}}P\).
  Here, since \(\induc_{G}^{\tilde{G}}P\) is a projective \(k\tilde{G}\)-module by \Cref{Theorem: Frobenius and projective} \ref{proj to proj 2022-06-02 12:44:53} and \(\Omega\tilde{M}\) has no projective summands by the self-injectivity of \(k\tilde{G}\),
  we have that \(\Omega\tilde{M}\) is a direct summand of \(\induc_{G}^{\tilde{G}}(\Omega(\restr_{G}^{\tilde{G}}\tilde{M}))\).
  Therefore, \(\restr_{G}^{\tilde{G}}(\Omega\tilde{M})\) is a direct summand of
  \begin{equation}
    \restr_{G}^{\tilde{G}}\induc_{G}^{\tilde{G}}(\Omega(\restr_{G}^{\tilde{G}}\tilde{M}))
    \cong \bigoplus_{\tilde{g}\in [\tilde{G}/G]} \tilde{g}\Omega(\restr_{G}^{\tilde{G}}\tilde{M})
  \end{equation}
  by \Cref{mackey}, which implies that \(\restr_{G}^{\tilde{G}}(\Omega\tilde{M})\) is has no projective summands because each \(\tilde{g}\Omega(\restr_{G}^{\tilde{G}}\tilde{M})\)
  has no projective summands by the self-injectivity of \(kG\).
  Thus, we conclude that \(P=0\) and
  \(\restr_{G}^{\tilde{G}}(\Omega\tilde{M})\cong \Omega(\restr_{G}^{\tilde{G}}\tilde{M})\).

  The later assertion follows from the fact that \(\tau \cong \Omega^2\)
  and the relative \(G\)-projectivity of \(\Omega\tilde{M}\).
\end{proof}

The following is important for the proof of \Cref{main-thm-1}.

\begin{proposition}[{\cite[Corollary 2.13]{MR3187626}}]\label{AIR Corollary2.13}
  Let \(\Lambda\) be a finite dimensional algebra.
  For a \(\tau\)-rigid pair \((M, P)\) for \(\Lambda\) the following are equivalent:
  \begin{enumerate}
    \item \((M, P)\) is a support \(\tau\)-tilting pair for \(\Lambda\).
    \item If \(\Hom_{\Lambda}(M, \tau X)=0, \Hom_{\Lambda}(X, \tau M)=0\) and \(\Hom_{\Lambda}(P, X)=0\), then \(X \in \add M\).
  \end{enumerate}
\end{proposition}

\begin{theorem}\label{main-thm-1}
  Let \(\tilde{M}\) be a support \(\tau\)-tilting \(k\tilde{G}\)-module.
  If it holds that \(\induc_{G}^{\tilde{G}}\restr_{G}^{\tilde{G}}\tilde{M} \in \add \tilde{M}\)
  and \(\tilde{M}\) is relatively \(G\)-projective,
  then we have that \(\restr_{G}^{\tilde{G}}\tilde{M}\) is a support \(\tau\)-tilting \(kG\)-module.
  Moreover, if \((\tilde{M}, \tilde{P})\) is a support \(\tau\)-tilting pair for \(k\tilde{G}\)
  corresponding to \(\tilde{M}\), then
  \((\restr_{G}^{\tilde{G}}\tilde{M}, \restr_{G}^{\tilde{G}}\tilde{P})\) is
  a support \(\tau\)-tilting pair for \(kG\)
  corresponding to \(\restr_{G}^{\tilde{G}}\tilde{M}\).
\end{theorem}

\begin{proof}
  Let \((\tilde{M}, \tilde{P})\) be a support \(\tau\)-tilting pair for \(k\tilde{G}\)
  corresponding to the support \(\tau\)-tilting \(k\tilde{G}\)-module \(\tilde{M}\).

  First, we show that \((\restr_{G}^{\tilde{G}}\tilde{M}, \restr_{G}^{\tilde{G}}\tilde{P})\)
  is a \(\tau\)-rigid pair for \(kG\).
  Since the \(k\tilde{G}\)-module \(\tilde{M}\) is a support \(\tau\)-tilting module,
  it is a \(\tau\)-rigid module.
  Hence, we have that \(\restr_{G}^{\tilde{G}}\tilde{M}\) is a \(\tau\)-rigid \(k\tilde{G}\)-module by \Cref{res-rigid}.
  On the other hand, by \Cref{Theorem: Frobenius and projective} we have that
  \begin{equation}
    \Hom_{kG}(\restr_{G}^{\tilde{G}}\tilde{P}, \restr_{G}^{\tilde{G}}\tilde{M})
    \cong \Hom_{k\tilde{G}}(\tilde{P}, \induc_{G}^{\tilde{G}}\restr_{G}^{\tilde{G}}\tilde{M}).
  \end{equation}
  Now, by the assumption that \(\induc_{G}^{\tilde{G}}\restr_{G}^{\tilde{G}}\tilde{M} \in \add \tilde{M}\),
  we have that \(\Hom_{kG}(\restr_{G}^{\tilde{G}}\tilde{P}, \restr_{G}^{\tilde{G}}\tilde{M})=0\)
  because \((\tilde{M}, \tilde{P})\) is a support \(\tau\)-tilting pair for \(k\tilde{G}\).
  Therefore, we conclude that
  \((\restr_{G}^{\tilde{G}}\tilde{M}, \restr_{G}^{\tilde{G}}\tilde{P})\)
  is a \(\tau\)-rigid pair for \(kG\).

  Next, we show that the \(\tau\)-rigid pair \((\restr_{G}^{\tilde{G}}\tilde{M}, \restr_{G}^{\tilde{G}}\tilde{P})\)
  is a support \(\tau\)-tilting pair for \(kG\).
  We show that \(X\in \add\restr_G^{\tilde{G}}\tilde{M}\)
  under the assumption that
  \begin{equation}
    \Hom_{kG}(X, \tau(\restr_G^{\tilde{G}}\tilde{M}))=\Hom_{kG}(\restr_G^{\tilde{G}}\tilde{M}, \tau X)=\Hom_{kG}(\restr_G^{\tilde{G}}\tilde{P}, X)=0,
  \end{equation}
  which implies that the pair \((\restr_{G}^{\tilde{G}}\tilde{M}, \restr_{G}^{\tilde{G}}\tilde{P})\)
  is a support \(\tau\)-tilting pair for \(kG\) by \Cref{AIR Corollary2.13}.
  Under these assumptions, we have the following:
  \begin{align}
    \Hom_{k\tilde{G}}(\induc_{G}^{\tilde{G}}X, \tau\tilde{M})
     & \cong \Hom_{kG}(X, \restr_{G}^{\tilde{G}}(\tau\tilde{M})) & \text{(\Cref{Theorem: Frobenius and projective})} \\
     & \cong \Hom_{kG}(X, \tau(\restr_{G}^{\tilde{G}}\tilde{M})) & \text{(\Cref{tau-restriction})}                   \\
     & =0.
  \end{align}
  \begin{align}
    \Hom_{k\tilde{G}}(\tilde{M}, \tau(\induc_{G}^{\tilde{G}}X))
     & \cong \Hom_{k\tilde{G}}(\tilde{M}, \induc_{G}^{\tilde{G}}(\tau X)) & \text{(\Cref{commute})}                           \\
     & \cong \Hom_{kG}(\restr_{G}^{\tilde{G}}\tilde{M}, \tau X)           & \text{(\Cref{Theorem: Frobenius and projective})} \\
     & =0.
  \end{align}
  \begin{align}
    \Hom_{k\tilde{G}}(\tilde{P}, \induc_{G}^{\tilde{G}}X)
     & \cong \Hom_{k\tilde{G}}(\restr_{G}^{\tilde{G}}\tilde{P}, X) & \text{(\Cref{Theorem: Frobenius and projective})} \\
     & =0.
  \end{align}

  By these three isomorphisms and the fact that \((\tilde{M}, \tilde{P})\)
  is a support \(\tau\)-tilting pair for \(k\tilde{G}\),
  applying \Cref{AIR Corollary2.13},
  we have that \( \induc_{G}^{\tilde{G}}X\in \add \tilde{M}\).
  Also, \(X\) is a direct summand of \( \restr_{G}^{\tilde{G}}\induc_{G}^{\tilde{G}}X\)
  by \Cref{mackey}.
  Therefore, we have that \(X\in \add \restr_{G}^{\tilde{G}}\tilde{M}\).
\end{proof}

\begin{corollary}\label{cor order 2022-12-26 18:27:13}
  Let \(\tilde{M}_1\) and \(\tilde{M}_2\) be relatively \(G\)-projective support \(\tau\)-tilting \(k\tilde{G}\)-modules such that \(\induc_{G}^{\tilde{G}}\restr_{G}^{\tilde{G}}\tilde{M}_i\in \add \tilde{M}_i\) for \(i=1,2\).
  Then \(\tilde{M}_1\geq \tilde{M}_2\) in \(\stautilt k\tilde{G}\) means that \(\restr_{G}^{\tilde{G}}\tilde{M}_1\geq \restr_{G}^{\tilde{G}}\tilde{M}_2\) in \(\stautilt kG\).
\end{corollary}
\begin{proof}
  The consequence immediately follows from \Cref{main-thm-1} and the exactness of the functor \(\restr_{G}^{\tilde{G}}\).
\end{proof}
We consider equivalent conditions to the assumption of \Cref{main-thm-1}.
First, we give the \namecrefs{dec lemma group algebra 2022-12-07 11:51:53} which can be applied in case of rigid \(k\tilde{G}\)-modules not only support \(\tau\)-tilting \(k\tilde{G}\)-modules.

\begin{lemma}\label{dec lemma group algebra 2022-12-07 11:51:53}
  Let \(\tilde{M}\) be a rigid \(k\tilde{G}\)-module and \(L\) a \(k\tilde{G}\)-module.
  If it holds that \(S\otimes_k \tilde{M}\in \add \tilde{M}\) for any composition factor \(S\) of \(L\),
  then the following isomorphism as \(k\tilde{G}\)-modules holds:
  \begin{equation}
    L\otimes_k \tilde{M} \cong \bigoplus_{S}S \otimes_k \tilde{M},
  \end{equation}
  where \(S\) runs over all composition factors of \(L\).
\end{lemma}
\begin{proof}
  Let \(L\) be an arbitrary \(k\tilde{G}\)-module and \(\tilde{M}\) a rigid \(k\tilde{G}\)-module satisfying that
  \begin{equation}
    S\otimes_k \tilde{M}\in \add \tilde{M} \text{ for any composition factor of \(S\) of \(L\)}.\label{dec lemma group algebra 2022-12-07 11:51:53 assumption}
  \end{equation}
  We use induction on the composition length \(\ell(L)\) of \(L\).
  If \(\ell(L)=1\), there is nothing to prove.
  Hence, we assume that \(\ell(L)\geq 2\) and that the statement for any \(k\tilde{G}\)-module \(L'\) satisfying \(\ell(L')<\ell(L)\) is true.
  Let \(T\) be a simple submodule of \(L\).
  We get the exact sequence
  \begin{equation}\label{ex_seq tensor 2022-12-07 12:48:07}
    \begin{tikzcd}
      0\ar[r]&T\otimes_k \tilde{M}\ar[r]&L\otimes_k \tilde{M}\ar[r]&L/T\otimes_k \tilde{M}\ar[r]&0
    \end{tikzcd}
  \end{equation}
  obtained by applying the exact functor \(-\otimes_k \tilde{M}\)
  to the exact sequence
  \begin{equation}
    \begin{tikzcd}
      0\ar[r]&T\ar[r]&L\ar[r]&L/T\ar[r]&0.
    \end{tikzcd}
  \end{equation}
  By the rigidity of \(\tilde{M}\), the assumption \eqref{dec lemma group algebra 2022-12-07 11:51:53 assumption} and the assumption of this induction, the sequence \eqref{ex_seq tensor 2022-12-07 12:48:07} splits, and we get that
  \begin{equation}
    L\otimes_k \tilde{M}\cong T\otimes_k \tilde{M}\oplus L/T\otimes_k \tilde{M} \cong T\otimes_k \tilde{M}\oplus \bigoplus_{S'}S'\otimes_k \tilde{M}\cong \bigoplus_{S}S\otimes_k \tilde{M},
  \end{equation}
  where \(S'\) and \(S\) run over all composition factors of \(L/T\) and \(L\), respectively.
\end{proof}

\begin{lemma}\label{equivalence-condition}
  Let \(\tilde{M}\) be a rigid \(k\tilde{G}\)-module.
  Then the following conditions are equivalent:
  \begin{enumerate}
    \item \(\induc_{G}^{\tilde{G}}\restr_{G}^{\tilde{G}} \tilde{M} \in \add \tilde{M}\) and \(\tilde{M}\) is relatively \(G\)-projective.\label{equivalence-condition item1}
    \item \(S\otimes_k \tilde{M} \in \add \tilde{M}\) for each simple \(k(\tilde{G}/G)\)-module \(S\).\label{equivalence-condition item2}
  \end{enumerate}
\end{lemma}
\begin{proof}
  By \Cref{Theorem: Frobenius and projective}, we have that
  \begin{equation}
    \induc_{G}^{\tilde{G}}\restr_{G}^{\tilde{G}} \tilde{M} \cong \induc_{G}^{\tilde{G}} (k_{G} \otimes_k \restr_{G}^{\tilde{G}}\tilde{M})
    \cong (\induc_{G}^{\tilde{G}} k_{G}) \otimes_k \tilde{M}
    \cong k(\tilde{G}/G) \otimes_k \tilde{M}.
  \end{equation}
  \ref{equivalence-condition item1} \(\Rightarrow\) \ref{equivalence-condition item2}.
  By the assumptions, we have that \(\induc_{G}^{\tilde{G}}\restr_{G}^{\tilde{G}} \tilde{M}=_{\add}\tilde{M}\).
  Hence, by \Cref{Theorem: Frobenius and projective} we get that
  \begin{align}
    S\otimes_k \tilde{M}
     & =_{\add}S\otimes_k \induc_{G}^{\tilde{G}}\restr_{G}^{\tilde{G}}\tilde{M}                        \\
     & \cong  \induc_{G}^{\tilde{G}}(\restr_{G}^{\tilde{G}}S\otimes_k\restr_{G}^{\tilde{G}}\tilde{M})  \\
     & \cong \induc_{G}^{\tilde{G}}( k_{G}^{\oplus \dim_k S} \otimes_k\restr_{G}^{\tilde{G}}\tilde{M}) \\
     & =_{\add} \induc_{G}^{\tilde{G}}\restr_{G}^{\tilde{G}} \tilde{M}                                 \\
     & =_{\add} \tilde{M},
  \end{align}
  for any simple \(k(\tilde{G}/G)\)-module \(S\),
  which implies that  \(S\otimes_k \tilde{M}\in \add \tilde{M}\).

  \noindent
  \ref{equivalence-condition item2} \(\Rightarrow\) \ref{equivalence-condition item1}.
  By \Cref{dec lemma group algebra 2022-12-07 11:51:53}, we have that
  \begin{equation}
    \induc_{G}^{\tilde{G}}\restr_{G}^{\tilde{G}} \tilde{M}\cong k(\tilde{G}/G) \otimes_k \tilde{M} \cong \bigoplus_{S}S \otimes_k \tilde{M},
  \end{equation}
  where \(S\) runs over all composition factors of the \(k\tilde{G}\)-module \(k(\tilde{G}/G)\).
  Therefore, the assumption implies that \(\induc_{G}^{\tilde{G}}\restr_{G}^{\tilde{G}} \tilde{M}\in \add \tilde{M}\).
  Moreover, since the trivial \(k\tilde{G}\)-module \(k_{\tilde{G}}\) appears as a composition factor of \(k(\tilde{G}/G)\), we have that the module \(\tilde{M}\) appears as a direct summand of \(\induc_{G}^{\tilde{G}}\restr_{G}^{\tilde{G}} \tilde{M}\), that is \(\tilde{M}\) is a relatively \(G\)-projective \(k\tilde{G}\)-module.
\end{proof}

We give the equivalent conditions to the assumption of \Cref{main-thm-1}.

\begin{theorem}\label{main-thm-2}
  Let \(\tilde{M}\) be a support \(\tau\)-tilting \(k\tilde{G}\)-module.
  Then the following conditions are equivalent:
  \begin{enumerate}
    \item \(\tilde{M}=_{\add} \induc_{G}^{\tilde{G}} M\) for some \(\tilde{G}\)-invariant support \(\tau\)-tilting \(kG\)-module \(M\).\label{main-thm-2 item1}
    \item \(\induc_{G}^{\tilde{G}}\restr_{G}^{\tilde{G}} \tilde{M} \in \add \tilde{M}\) and \(\tilde{M}\) is relatively \(G\)-projective.\label{main-thm-2 item2}
    \item \(S\otimes_k \tilde{M} \in \add \tilde{M}\) for each simple \(k(\tilde{G}/G)\)-module \(S\).\label{main-thm-2 item3}
  \end{enumerate}
\end{theorem}
\begin{proof}
  \ref{main-thm-2 item1} \(\Rightarrow\) \ref{main-thm-2 item2}.
  Assume that \(\tilde{M}=_{\add} \induc_{G}^{\tilde{G}} M\) for some \(\tilde{G}\)-invariant
  support \(\tau\)-tilting \(kG\)-module \(M\).
  Then clearly \(\tilde{M}\) is a relatively \(G\)-projective \(k\tilde{G}\)-module (see \cite[3.9.1]{MR860771}), and by \Cref{mackey}, we have that
  \begin{equation}
    \induc_{G}^{\tilde{G}}\restr_{G}^{\tilde{G}} \tilde{M}=_{\add}\induc_{G}^{\tilde{G}}\restr_{G}^{\tilde{G}} \induc_{G}^{\tilde{G}} M \cong \induc_{G}^{\tilde{G}}(\bigoplus_{\tilde{g}\in [\tilde{G}/G]}\tilde{g}M)
    \cong \bigoplus_{\tilde{g}\in [\tilde{G}/G]}\induc_{G}^{\tilde{G}} M \in \add \tilde{M}.
  \end{equation}

  \noindent
  \ref{main-thm-2 item2} \(\Rightarrow\) \ref{main-thm-2 item1}.
  Assume that \(\induc_{G}^{\tilde{G}}\restr_{G}^{\tilde{G}}\tilde{M} \in \add \tilde{M}\)
  and that \(\tilde{M}\) is relatively \(G\)-projective.
  Put \(M:=\restr_{G}^{\tilde{G}}\tilde{M}\).
  Then by \Cref{restriction-invariance}
  and \Cref{main-thm-1},
  \(M\) is a \(\tilde{G}\)-invariant support \(\tau\)-tilting \(kG\)-module.
  We show that \(\induc_{G}^{\tilde{G}} M =_{\add}\tilde{M}\), that is \(\add(\induc_{G}^{\tilde{G}} M) =\add\tilde{M}\).
  By the assumption that \(\induc_{G}^{\tilde{G}}\restr_{G}^{\tilde{G}}\tilde{M}\in \add \tilde{M}\),
  we have \(\add(\induc_{G}^{\tilde{G}} M) \subset\add\tilde{M}\).
  On the other hand, since \(\tilde{M}\) is relatively \(G\)-projective,
  \(\tilde{M}\) is a direct summand of \(\induc_{G}^{\tilde{G}}\restr_{G}^{\tilde{G}} \tilde{M}=\induc_{G}^{\tilde{G}} M\).
  Hence, we have \(\add\tilde{M} \subset \add(\induc_{G}^{\tilde{G}} M)\).

  \noindent
  \ref{main-thm-2 item2} \(\Leftrightarrow\) \ref{main-thm-2 item3}.
  Since support \(\tau\)-tilting \(k\tilde{G}\)-modules
  are rigid \(k\tilde{G}\)-modules,
  the equivalence
  follows from \Cref{equivalence-condition}.
\end{proof}

\begin{corollary}\label{cor1}
  Let \((\stautilt kG)^{\tilde{G}}\) be the subset of
  \(\stautilt kG\) consisting of \(\tilde{G}\)-invariant support \(\tau\)-tilting
  \(kG\)-modules
  and
  \((\stautilt k\tilde{G})^{\star}\) the subset of
  \(\stautilt k\tilde{G}\) consisting of support \(\tau\)-tilting
  \(k\tilde{G}\)-modules satisfying the equivalent conditions of Theorem
  \ref{main-thm-2}.
  Then the induction functor \(\induc_{G}^{\tilde{G}}\) induces
  a poset isomorphism
  \begin{equation}\label{bij 2022-12-13 14:06:34}
    \begin{tikzcd}[row sep=0pt]
      (\stautilt kG)^{\tilde{G}}\ar[r,"\sim"]& (\stautilt k\tilde{G})^{\star}\\
      M\ar[r,mapsto]&\induc_{G}^{\tilde{G}} M.
    \end{tikzcd}
  \end{equation}
  In particular, the induction functor \(\induc_{G}^{\tilde{G}}\) induces the poset monomorphism
  \begin{equation}\label{mono 2022-12-13 14:04:30}
    \begin{tikzcd}[row sep=0pt]
      (\stautilt kG)^{\tilde{G}}\ar[r]& \stautilt k\tilde{G}\\
      M\ar[r,mapsto]&\induc_{G}^{\tilde{G}} M.
    \end{tikzcd}
  \end{equation}
\end{corollary}

\begin{proof}
  By  \cite[Theorem 3.2]{https://doi.org/10.48550/arxiv.2208.14680},
  the map \eqref{mono 2022-12-13 14:04:30} is well-defined.
  % Also, the map is a poset homomorphism,
  Moreover, by the exactness of the functor \(\induc_{G}^{\tilde{G}}\), if \(N\leq M\) in \(\stautilt kG\) then \(\induc_{G}^{\tilde{G}} N \leq \induc_{G}^{\tilde{G}} M\) in \(\stautilt k\tilde{G}\) for any support \(\tau\)-tilting \(kG\)-modules \(N\) and \(M\).
  Therefore, the map \eqref{mono 2022-12-13 14:04:30} is a poset homomorphism.
  % because for a epimorphism
  % \(
  % \begin{tikzcd}[column sep=20pt]
  %   M^{\oplus r}\ar[r,->>]& N
  % \end{tikzcd}
  % \),
  % the induced homomorphism
  % \(
  % \begin{tikzcd}[column sep=20pt]
  %   (\induc_{G}^{\tilde{G}} M)^{\oplus r} \ar[r,->>]& \induc_{G}^{\tilde{G}} N
  % \end{tikzcd}
  % \)
  % is an epimorphism,
  % which means that if \(N\leq M\) then \(\induc_{G}^{\tilde{G}} N \leq \induc_{G}^{\tilde{G}} M\) for any support \(\tau\)-tilting \(kG\)-modules \(N\) and \(M\).

  We show that the map \eqref{mono 2022-12-13 14:04:30} restricts to a bijection \eqref{bij 2022-12-13 14:06:34}.
  By the definition of \((\stautilt k\tilde{G})^{\star}\) and the above argument, the map \eqref{bij 2022-12-13 14:06:34} is well-defined and a poset homomorphism.
  For any relatively \(G\)-projective
  support \(\tau\)-tilting \(k\tilde{G}\)-module \(\tilde{M}\)
  with \(\induc_{G}^{\tilde{G}}\restr_{G}^{\tilde{G}}\tilde{M}\in \tilde{M}\),
  by \Cref{main-thm-2}, we can take a \(\tilde{G}\)-invariant support \(\tau\)-tilting \(kG\)-module \(M\) satisfying \(\induc_{G}^{\tilde{G}}M =_{\add}\tilde{M}\).
  Hence, the map is surjective.
  Also, assume that two \(\tilde{G}\)-invariant support \(\tau\)-tilting \(kG\)-modules
  \(M\) and \(N\) satisfy that \(\induc_{G}^{\tilde{G}} M =_{\add} \induc_{G}^{\tilde{G}} N\).
  Then we have that
  \begin{equation}
    \bigoplus_{\tilde{g}\in [\tilde{G}/G]} \tilde{g}M
    \cong \restr_{G}^{\tilde{G}}\induc_{G}^{\tilde{G}} M
      =_{\add} \restr_{G}^{\tilde{G}}\induc_{G}^{\tilde{G}} N
    \cong \bigoplus_{\tilde{g}\in [\tilde{G}/G]}\tilde{g}N
  \end{equation}
  by \Cref{mackey},
  which means that \(M=_{\add} N\) by the \(\tilde{G}\)-invariances of \(M\) and \(N\).
  Hence, the map is injective.
  This completes the proof of the first assertion.

  The latter assertion immediately follows from the fact that the map \eqref{mono 2022-12-13 14:04:30} is
  the composition of the poset isomorphism \eqref{bij 2022-12-13 14:06:34} and the inclusion map
  \(
  \begin{tikzcd}
    (\stautilt k\tilde{G})^{\star}\ar[r,hook]&\stautilt k\tilde{G}.
  \end{tikzcd}
  \)
\end{proof}

As an application of \Cref{main-thm-2},
we consider the case that \(\tilde{G}/G\) is a \(p\)-group.

\begin{theorem}\label{p-extension}
  Let \(\tilde{G}\) be a finite group and
  \(G\) a normal subgroup of \(\tilde{G}\)
  of \(p\)-power index in \(\tilde{G}\).
  Then the induction functor
  \(\induc_{G}^{\tilde{G}}\) induces an isomorphism as partially ordered sets
  between \((\stautilt kG)^{\tilde{G}}\) and \(\stautilt k\tilde{G}\),
  where \((\stautilt kG)^{\tilde{G}}\) is
  the subset of \(\stautilt kG\) consisting of \(\tilde{G}\)-invariant support \(\tau\)-tilting
  \(kG\)-module.
\end{theorem}

\begin{proof}
  By \Cref{cor1}, the map \eqref{bij 2022-12-13 14:06:34} is a poset isomorphism.
  We enough to show that \((\stautilt k\tilde{G})^{\star}=\stautilt k\tilde{G}\).
  It is clear that \((\stautilt k\tilde{G})^{\star}\subset\stautilt k\tilde{G}\).
  To prove the reverse inclusion,
  take an arbitrary support \(\tau\)-tilting \(k\tilde{G}\)-module \(\tilde{M}\).
  Since \(\tilde{G}/G\) is a \(p\)-group,
  the only simple \(k(\tilde{G}/G)\)-module is the trivial \(k(\tilde{G}/G)\)-module.
  Hence, the condition \ref{main-thm-2 item3} of \Cref{main-thm-2}
  is satisfied in our situation
  because \(k_{\tilde{G}/G}\otimes_k \tilde{M}\)
  is isomorphic to \(\tilde{M}\).
  Therefore, we conclude that
  \(\stautilt k\tilde{G}\subset (\stautilt k\tilde{G})^{\star}\).
\end{proof}
The following example can be seen in \cite{https://doi.org/10.48550/arxiv.2208.14680}.
\begin{example}\label{A4S4Example 2023-01-04 20:09:08}
  Let \(k\) be an algebraically closed field of characteristic \(p=2\).
  We consider that the case that \(G\) is the alternating group \(A_4\) of degree \(4\) and \(\tilde{G}\)
  is the symmetric group \(S_4\) of degree \(4\).
  The algebras \(kA_4\) and \(kS_4\) are Brauer graph algebras associated to the Brauer graphs in \Cref{A4 Brauer graph 2022-07-16 22:29:08} and \Cref{S4 Brauer graph 2022-07-16 22:28:53}, respectively:

  \begin{figure}[ht]
    \begin{center}
      \subfigure[The Brauer graph of \(kA_4\)]{
        \begin{tikzpicture}
          % \draw[help lines] (0,0) grid (5,3);
          % \coordinate (x) at (-10,0);
          \coordinate (O) at (0,0);
          \coordinate (A) at (4,0);
          \coordinate (B) at (2,3);
          \draw (O) edge node[below] {\(k_{A_4}=1\)} (A) ;
          \draw (A) edge node[right] {\(2\)} (B);
          \draw (B) edge node[left] {\(3\)}(O);
          \foreach \P in {O,A,B} \fill[white] (\P) circle (2pt);
          \foreach \P in {O,A,B} \draw (\P) circle (2pt);
        \end{tikzpicture}
        \label{A4 Brauer graph 2022-07-16 22:29:08}
      }
      \hfill
      \subfigure[The Brauer graph of \(kS_4\)]{

        \begin{tikzpicture}
          % \node(O) at (0,0){\(O\)};
          % \draw[help lines] (0,0) grid (10,3);
          \node(m) at (5.5,1.5){multiplicity: \(2\)};
          \coordinate (O) at (2,1.5);
          \coordinate (A) at (4,1.5);
          \coordinate (B) at (1,1.5);
          \draw (O) edge node[below] {\(2'\)} (A) ;
          \foreach \P in {O,A} \draw (\P) circle (2pt);
          \draw[out=-90,in=-90] (A) edge node[below] {\(1'=k_{S_4}\)} (B);
          \draw[out=90,in=90] (A) edge node[below] {} (B);
          \fill[white] (A) circle (2pt);
          \fill[black] (O) circle (2pt);
        \end{tikzpicture}

        \label{S4 Brauer graph 2022-07-16 22:28:53}
      }
    \end{center}
    \vspace{-0.5cm}
    \caption{Brauer graphs}
    \label{Brauer 2022-07-16 12:36:42}
  \end{figure}
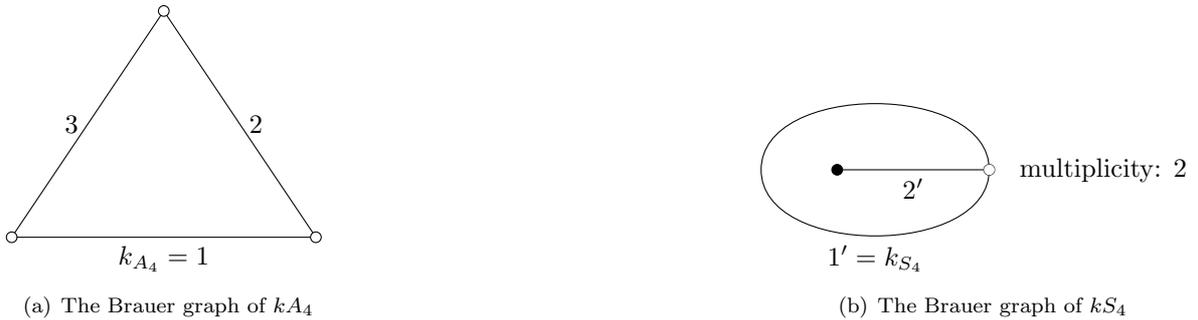

  Now we draw the Hasse diagram \(\mathcal{H}(\stautilt kA_4)\) of the partially ordered set \(\stautilt kA_4\) as follows:
  \begin{landscape}
    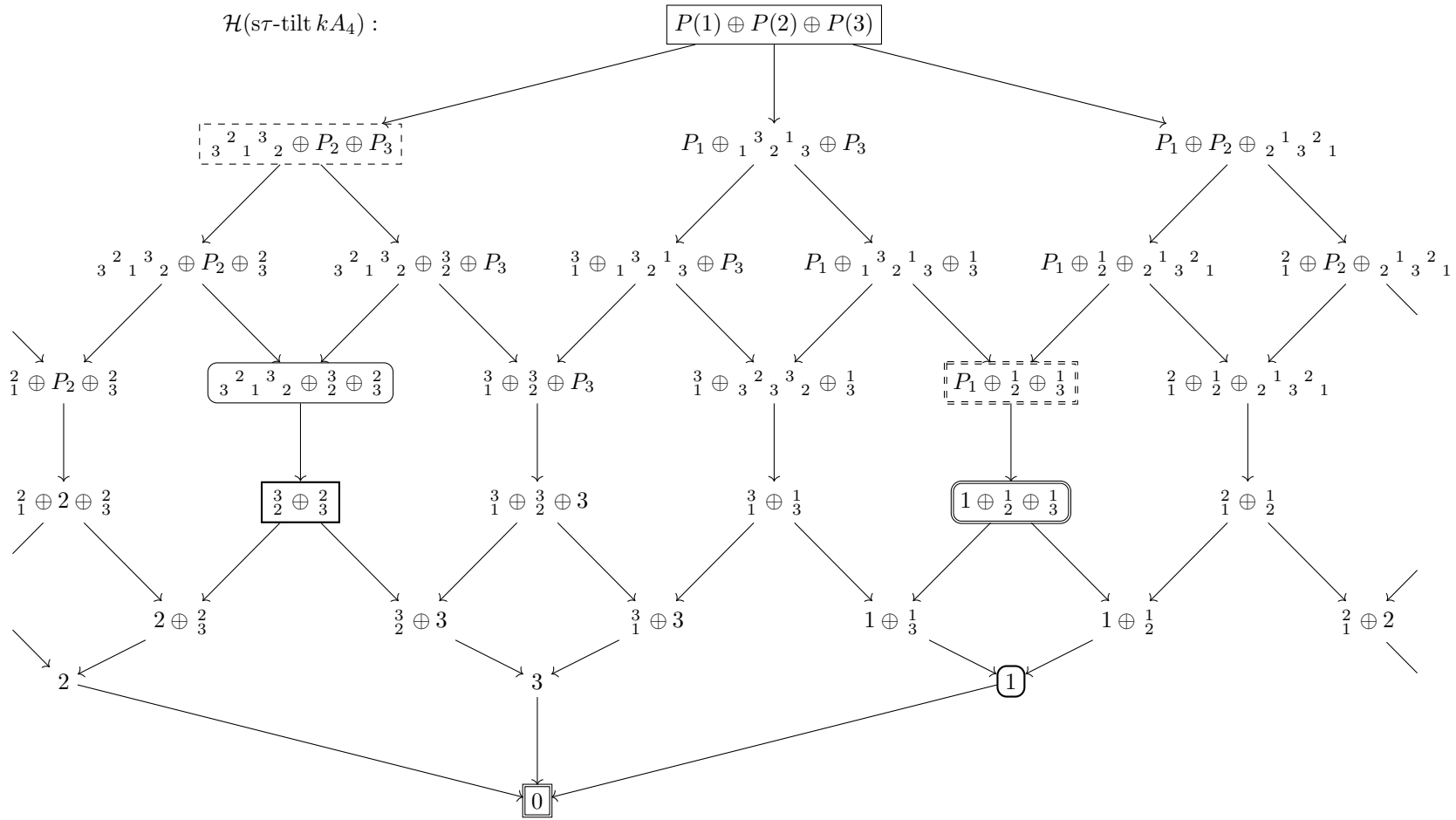
\begin{figure}[ht]
      \begin{tikzpicture}[scale=0.9]
        % \node(O) at (0,0){\(O\)};
        % \draw[help lines] (-10,-10) grid (15,6);
        \node(hasse) at (-4,6) {\(\mathcal{H}(\stautilt kA_4):\)};
        \node(progenerator)[draw] at (4,6){\(
          P(1)\oplus P(2)\oplus P(3)
          \)};
        \node(a) at (4,4){\(P_1\oplus
          \begin{smallmatrix}
            & 3 &   & 1 &   \\
            1 &   & 2 &   & 3 \\
          \end{smallmatrix}
          \oplus P_3\)};
        \node(j) at (4,0){\(
          \begin{smallmatrix}
            3\\
            1
          \end{smallmatrix}\oplus
          \begin{smallmatrix}
            & 2 &   & 3 &   \\
            3 &   & 3 &   & 2 \\
          \end{smallmatrix}\oplus
          \begin{smallmatrix}
            1\\
            3
          \end{smallmatrix}
          \)};
        \node(n) at (-0,0){\(
          \begin{smallmatrix}
            3\\
            1
          \end{smallmatrix}\oplus
          \begin{smallmatrix}
            3\\
            2
          \end{smallmatrix}\oplus
          P_3
          \)};
        \node(w) at (0,-2){\(
          \begin{smallmatrix}
            3\\
            1
          \end{smallmatrix}\oplus
          \begin{smallmatrix}
            3\\
            2
          \end{smallmatrix}\oplus
          3
          \)};
        \node(x) at (2,-4){\(
          \begin{smallmatrix}
            3\\
            1
          \end{smallmatrix}\oplus
          3\)};
        \node(z) at (-2,-4){\(
          \begin{smallmatrix}
            3\\
            2
          \end{smallmatrix}\oplus
          3
          \)};
        \node(y) at (0,-5){\(
          3
          \)};
        \node(k) at (4,-2){\(
          \begin{smallmatrix}
            3\\
            1
          \end{smallmatrix}\oplus
          \begin{smallmatrix}
            1\\
            3
          \end{smallmatrix}
          \)};
        \node(g) at (2,2){\(
          \begin{smallmatrix}
            3\\
            1
          \end{smallmatrix}\oplus
          \begin{smallmatrix}
            & 3 &   & 1 &   \\
            1 &   & 2 &   & 3 \\
          \end{smallmatrix}
          \oplus P_3
          \)};
        \node(h) at (6,2){\(
          P_1\oplus
          \begin{smallmatrix}
            & 3 &   & 1 &   \\
            1 &   & 2 &   & 3 \\
          \end{smallmatrix}\oplus
          \begin{smallmatrix}
            1\\
            3
          \end{smallmatrix}
          \)};
        \node(b)[draw,dashed] at (-4,4){\(
          \begin{smallmatrix}
            & 2 &   & 3 &   \\
            3 &   & 1 &   & 2 \\
          \end{smallmatrix}
          \oplus P_2\oplus P_3
          \)};
        \node(c) at (12,4){\(
          P_1\oplus P_2\oplus
          \begin{smallmatrix}
            & 1 &   & 2 &   \\
            2 &   & 3 &   & 1 \\
          \end{smallmatrix}
          \)};
        \node(d) at (-6,2){\(
          \begin{smallmatrix}
            & 2 &   & 3 &   \\
            3 &   & 1 &   & 2 \\
          \end{smallmatrix}
          \oplus P_2\oplus
          \begin{smallmatrix}
            2\\3
          \end{smallmatrix}
          \)};
        \node(l)[draw,rounded corners] at (-4,0){\(
          \begin{smallmatrix}
            & 2 &   & 3 &   \\
            3 &   & 1 &   & 2 \\
          \end{smallmatrix}\oplus
          \begin{smallmatrix}
            3\\
            2
          \end{smallmatrix}\oplus
          \begin{smallmatrix}
            2\\
            3
          \end{smallmatrix}
          \)};
        \node(a1)[draw,thick] at (-4,-2){\(
          \begin{smallmatrix}
            3\\
            2
          \end{smallmatrix}\oplus
          \begin{smallmatrix}
            2\\
            3
          \end{smallmatrix}
          \)};
        \node(b1) at (-6,-4){\(
          2\oplus
          \begin{smallmatrix}
            2\\
            3
          \end{smallmatrix}
          \)};
        \node(c1) at (-8,-5){\(
          2
          \)};
        \node(m) at (-8,0){\(
          \begin{smallmatrix}
            2\\
            1
          \end{smallmatrix}\oplus
          P_2\oplus
          \begin{smallmatrix}
            2\\
            3
          \end{smallmatrix}
          \)};
        \node(pa1) at (-9,1){};
        \node(d1) at (-8,-2){\(
          \begin{smallmatrix}
            2\\
            1
          \end{smallmatrix}\oplus
          2\oplus
          \begin{smallmatrix}
            2\\
            3
          \end{smallmatrix}
          \)};
        \node(pa5) at (-9,-3){};
        \node(pa6) at (-9,-4){};
        \node(f) at (-2,2){\(
          \begin{smallmatrix}
            & 2 &   & 3 &   \\
            3 &   & 1 &   & 2 \\
          \end{smallmatrix}
          \oplus
          \begin{smallmatrix}
            3\\
            2
          \end{smallmatrix}
          \oplus
          P_3
          \)};
        \node(e) at (10,2){\(
          P_1\oplus
          \begin{smallmatrix}
            1\\
            2
          \end{smallmatrix}\oplus
          \begin{smallmatrix}
            & 1 &   & 2 &   \\
            2 &   & 3 &   & 1 \\
          \end{smallmatrix}
          \)};
        \node(u)[draw,double,dashed] at (8,0){\(
          P_1\oplus
          \begin{smallmatrix}
            1\\
            2
          \end{smallmatrix}\oplus
          \begin{smallmatrix}
            1\\
            3
          \end{smallmatrix}
          \)};
        \node(v)[draw,double,rounded corners] at (8,-2){\(
          1\oplus
          \begin{smallmatrix}
            1\\
            2
          \end{smallmatrix}\oplus
          \begin{smallmatrix}
            1\\
            3
          \end{smallmatrix}
          \)};
        \node(o) at (12,0){\(
          \begin{smallmatrix}
            2\\
            1
          \end{smallmatrix}\oplus
          \begin{smallmatrix}
            1\\
            2
          \end{smallmatrix}\oplus
          \begin{smallmatrix}
            & 1 &   & 2 &   \\
            2 &   & 3 &   & 1 \\
          \end{smallmatrix}
          \)};
        \node(p) at (12,-2){\(
          \begin{smallmatrix}
            2\\
            1
          \end{smallmatrix}\oplus
          \begin{smallmatrix}
            1\\
            2
          \end{smallmatrix}
          \)};
        \node(t) at (14,-4){\(
          \begin{smallmatrix}
            2\\
            1
          \end{smallmatrix}\oplus
          2
          \)};
        \node(pa3) at (15,-3){};
        \node(pa4) at (15,-5){};
        \node(q) at (10,-4){\(
          1
          \oplus
          \begin{smallmatrix}
            1\\
            2
          \end{smallmatrix}
          \)};
        \node(r) at (6,-4){\(
          1
          \oplus
          \begin{smallmatrix}
            1\\
            3
          \end{smallmatrix}
          \)};
        \node(s)[draw,thick,rounded corners] at (8,-5){\(
          1
          \)};
        \node(i) at (14,2){\(
          \begin{smallmatrix}
            2\\
            1
          \end{smallmatrix}\oplus
          P_2\oplus
          \begin{smallmatrix}
            & 1 &   & 2 &   \\
            2 &   & 3 &   & 1 \\
          \end{smallmatrix}
          \)};
        \node(pa2) at (15,1){};
        \node(zero)[draw,double] at (0,-7){\(0\)};
        \draw[->](progenerator) edge node{} (a);
        \draw[->](progenerator) edge node{} (b);
        \draw[->](progenerator) edge node{} (c);
        \draw[->](b) edge node{} (d);
        \draw[->](a) edge node{} (g);
        \draw[->](a) edge node{} (h);
        \draw[->](b) edge node{} (f);
        \draw[->](d) edge node{} (l);
        \draw[->](d) edge node{} (m);
        \draw[->](f) edge node{} (l);
        \draw[->](c) edge node{} (e);
        \draw[->](c) edge node{} (i);
        \draw[->](y) edge node{} (zero);
        \draw[->](p) edge node{} (q);
        \draw[->](q) edge node{} (s);
        \draw[->](r) edge node{} (s);
        \draw[->](s) edge node{} (zero);
        \draw[->](c1) edge node{} (zero);
        \draw[->](b1) edge node{} (c1);
        \draw[->](i) edge node{} (o);
        \draw[->](e) edge node{} (o);
        \draw[->](e) edge node{} (u);
        \draw[->](o) edge node{} (p);
        \draw[->](x) edge node{} (y);
        \draw[->](z) edge node{} (y);
        \draw[->](g) edge node{} (n);
        \draw[->](f) edge node{} (n);
        \draw[->](h) edge node{} (j);
        \draw[->](h) edge node{} (u);
        \draw[->](u) edge node{} (v);
        \draw[->](v) edge node{} (q);
        \draw[->](v) edge node{} (r);
        \draw[->](k) edge node{} (r);
        \draw[->](k) edge node{} (x);
        \draw[->](g) edge node{} (j);
        % \draw[->](i) edge node{} (m);
        \draw[->](m) edge node{} (d1);
        \draw[->](l) edge node{} (a1);
        \draw[->](w) edge node{} (x);
        \draw[->](w) edge node{} (z);
        \draw[->](a1) edge node{} (z);
        \draw[->](a1) edge node{} (b1);
        \draw[->](d1) edge node{} (b1);
        \draw[->](p) edge node{} (t);
        \draw[->](n) edge node{} (w);
        \draw[->](j) edge node{} (k);
        \draw[->](pa1) edge node{} (m);
        \draw[-](i) edge node{} (pa2);
        \draw[->](pa3) edge node{} (t);
        \draw[-](t) edge node{} (pa4);
        \draw[-](d1) edge node{} (pa5);
        \draw[->](pa6) edge node{} (c1);
      \end{tikzpicture}
      \caption[]{The Hasse diagram of \(\stautilt kA_4\)}
      \label{Hasse kA4 2022-07-07 02:43:05}
    \end{figure}
  \end{landscape}
  \noindent
  The enclosed support \(\tau\)-tilting modules in \Cref{Hasse kA4 2022-07-07 02:43:05} are all the invariant support \(\tau\)-tilting modules under the action of \(S_4\).
  Next, we draw the Hasse diagram \(\mathcal{H}(\stautilt kS_4)\) of partially ordered set \(\stautilt kS_4\) as follows:
  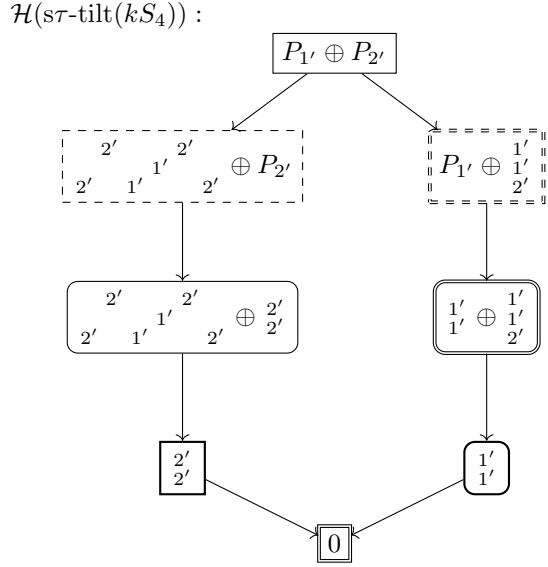
\begin{figure}[ht]
    \centering
    \begin{tikzpicture}
      \node(H) at (-5,1){\(\mathcal{H}(\stautilt(kS_4)):\)};
      \node(a)[draw] at (-2,0.5){\(
        P_{1'}\oplus P_{2'}
        \)};
      \node(b)[draw,dashed] at (-4,-1){\(
        \begin{smallmatrix}
          &2'&&&2'&\\
          &&&1'&&\\
          2'&&1'&&&2'
        \end{smallmatrix}
        \oplus P_{2'}
        \)};
      \node(d)[draw,rounded corners] at (-4,-3){\(
        \begin{smallmatrix}
          &2'&&&2'&\\
          &&&1'&&\\
          2'&&1'&&&2'
        \end{smallmatrix}
        \oplus
        \begin{smallmatrix}
          2'\\
          2'
        \end{smallmatrix}
        \)};
      % \node(e)[fill=olive] at (-4,-5){\(
      \node(e)[draw,thick] at (-4,-5){\(
        \begin{smallmatrix}
          2'\\
          2'
        \end{smallmatrix}
        \)};
      \node(zero)[draw,double] at (-2,-6){\(
        0
        \)};
      \node(c)[draw,double,dashed] at (0,-1){\(
        P_{1'}\oplus
        \begin{smallmatrix}
          1'\\
          1'\\
          2'
        \end{smallmatrix}
        \)};
      \node(f)[draw,double,rounded corners] at (0,-3){\(
        \begin{smallmatrix}
          1'\\
          1'
        \end{smallmatrix}
        \oplus
        \begin{smallmatrix}
          1'\\
          1'\\
          2'
        \end{smallmatrix}
        \)};
      \node(g)[draw,thick,rounded corners] at (0,-5){\(
        \begin{smallmatrix}
          1'\\
          1'
        \end{smallmatrix}
        \)};
      \draw[->](a) edge node{} (b);
      \draw[->](a) edge node{} (c);
      \draw[->](b) edge node{} (d);
      \draw[->](d) edge node{} (e);
      \draw[->](e) edge node{} (zero);
      \draw[->](g) edge node{} (zero);
      \draw[->](f) edge node{} (g);
      \draw[->](c) edge node{} (f);
    \end{tikzpicture}
    \caption{The Hasse diagram of \(\stautilt kS_4\)}
    \label{Hasse kS4 2022-07-07 02:39:33}
  \end{figure}

  \noindent
  The functor \(\induc_{A_4}^{S_4}\) takes each enclosed \(S_4\)-invariant support \(\tau\)-tilting \(kA_4\)-module in \Cref{Hasse kA4 2022-07-07 02:43:05} to that in \Cref{Hasse kS4 2022-07-07 02:39:33} with the same square.
\end{example}
\begin{remark}
  Let $(\induc_{G}^{\tilde{G}})^{-1}(\stautilt k\tilde{G})
    :=\{ M\in \stautilt kG\:|\: \induc_{G}^{\tilde{G}}M\in \stautilt k\tilde{G} \}$.
  Then \((\stautilt kG)^{\tilde{G}}\) is contained in \((\induc_{G}^{\tilde{G}})^{-1}(\stautilt k\tilde{G})\)
  by \cite[Theorem 3.2]{https://doi.org/10.48550/arxiv.2208.14680}.
  On the other hand, they do not coincide in general.
  Moreover, though the poset homomorphism
  \begin{equation}
    \begin{tikzcd}[row sep = 0pt]
      (\stautilt kG)^{\tilde{G}} \ar[r]& \stautilt k\tilde{G}\\
      M\ar[r,mapsto]&\induc_{G}^{\tilde{G}}M
    \end{tikzcd}
  \end{equation}
  % \(\induc_{G}^{\tilde{G}}: (\stautilt kG)^{\tilde{G}} \rightarrow \stautilt k\tilde{G}\)
  is a monomorphism by \Cref{cor1},
  the one
  \begin{equation}
    \begin{tikzcd}[row sep = 0pt]
      (\induc_{G}^{\tilde{G}})^{-1}(\stautilt k\tilde{G}) \ar[r]& \stautilt k\tilde{G}\\
      M\ar[r,mapsto]&\induc_{G}^{\tilde{G}}M
    \end{tikzcd}
  \end{equation}
  % \(\induc_{G}^{\tilde{G}}: (\induc_{G}^{\tilde{G}})^{-1}(\stautilt k\tilde{G}) \rightarrow \stautilt k\tilde{G}\) 
  is not a monomorphism in general.

  For example, for \(p=2\), the alternating group \(A_4\) of degree \(4\) and the symmetric group \(S_4\) of degree \(4\), a \(kA_4\)-module
  \( M:=
  1
  \oplus
  \begin{smallmatrix}
    1\\
    2
  \end{smallmatrix}
  \)
  is a support \(\tau\)-tilting \(kA_4\)-module, where \(1\) means the trivial \(kA_4\)-module and \(2\) a non-trivial \(kA_4\)-module.
  Also, it holds that
  \( \sigma M\cong
  1
  \oplus
  \begin{smallmatrix}
    1\\
    3
  \end{smallmatrix}
  \)
  for \(\sigma \in S_4\setminus A_4\), where \(3\) means the non-trivial simple \(kA_4\)-module not isomorphic to \(2\).
  Therefore, we have that \(M \not\in (\stautilt kA_4)^{S_4}\).
  However,
  \(
  \induc_{A_4}^{S_4}M \cong
  \begin{smallmatrix}
    1'\\
    1'
  \end{smallmatrix}          \oplus
  \begin{smallmatrix}
    1'\\
    1'\\
    2'
  \end{smallmatrix}
  \)
  is a support \(\tau\)-tilting $kS_4$-module, where \(1'\) means the trivial \(kS_4\)-module and \(2'\) the simple \(kS_4\)-module of dimension 2.
  This implies that \(M \in  (\induc_{A_4}^{S_4})^{-1}(\stautilt kS_4)\).
  Moreover, for
  \(
  N:=
  1
  \oplus
  \begin{smallmatrix}
    1\\
    2
  \end{smallmatrix}
  \oplus
  \begin{smallmatrix}
    1\\
    3
  \end{smallmatrix}
  \in (\stautilt kA_4)^{S_4}
  \),
  it holds that
  \(
  \induc_{A_4}^{S_4}N
  \cong
  \begin{smallmatrix}
    1'\\
    1'
  \end{smallmatrix}\oplus
  \begin{smallmatrix}
    1'\\
    1'\\
    2'
  \end{smallmatrix}\oplus
  \begin{smallmatrix}
    1'\\
    1'\\
    2'
  \end{smallmatrix}
  =_{\add}
  \begin{smallmatrix}
    1'\\
    1'
  \end{smallmatrix}          \oplus
  \begin{smallmatrix}
    1'\\
    1'\\
    2'
  \end{smallmatrix}
  (=_{\add}\induc_{A_4}^{S_4}M)
  \).
  Therefore, the map
  \begin{equation}
    \begin{tikzcd}[row sep = 0pt]
      (\induc_{A_4}^{S_4})^{-1}(\stautilt kS_4) \ar[r]& \stautilt kS_4\\
      M\ar[r,mapsto]&\induc_{A_4}^{S_4}M
    \end{tikzcd}
  \end{equation}
  is not a monomorphism.
\end{remark}

At the end of this section, we discuss a feature of vertices of indecomposable \(\tau\)-rigid \(k\tilde{G}\)-modules.

\begin{lemma}\label{Landrock}
  Let \(\tilde{G}\) be a finite group.
  Then the trivial \(k\tilde{G}\)-module \(k_{\tilde{G}}\) is a \(\tau\)-rigid if and only if \(\tilde{G}\) has no normal subgroup of index \(p\).
\end{lemma}

\begin{proof}
  By \cite[Chap. I, Corollary 10.13]{MR737910},
  there exists a normal subgroup of \(\tilde{G}\) of index \(p\)
  if and only if \(\Ext_{k\tilde{G}}^1(k_{\tilde{G}}, k_{\tilde{G}}) \neq 0\).
  Also, by the simplicity of the trivial \(k\tilde{G}\)-module and  Auslander-Reiten duality,
  we have that
  \begin{equation}
    \Hom_{k\tilde{G}}(k_{\tilde{G}}, \tau k_{\tilde{G}})
    \cong
    \overline{\Hom}_{k\tilde{G}}(k_{\tilde{G}}, \tau k_{\tilde{G}})
    \cong
    D\Ext_{k\tilde{G}}^1(k_{\tilde{G}}, k_{\tilde{G}}).
  \end{equation}
  Therefore, we get the result.
\end{proof}

\begin{theorem}\label{tau-rigid}
  Let \(\tilde{G}\) be a finite group.
  Then any indecomposable \(\tau\)-rigid \(k\tilde{G}\)-module has a vertex contained in
  a Sylow \(p\)-subgroup of \(\tilde{G}\) properly
  if and only if \(\tilde{G}\) has a proper normal subgroup of \(p\)-power index.
\end{theorem}

\begin{proof}
  Assume that \(\tilde{G}\) has no proper normal subgroup of \(p\)-power index.
  Then by \Cref{Landrock}, the trivial \(k\tilde{G}\)-module,
  whose vertex is a Sylow \(p\)-subgroup of \(\tilde{G}\), is a \(\tau\)-rigid module.

  Conversely, assume that \(\tilde{G}\) has normal subgroup of \(p\)-power index.
  In this case, there exists a normal subgroup \(G\) of \(\tilde{G}\) of index \(p\).
  Let \(\tilde{X}\) be an arbitrary \(\tau\)-rigid \(k\tilde{G}\)-module.
  Then, \(\tilde{X}\) is a direct summand of a support \(\tau\)-tilting \(k\tilde{G}\)-module \(\tilde{M}\) by \cite[Theorem 2.10]{MR3187626}, that is, \(\tilde{X}\) is relatively \(G\)-projective.
  Also, there exists a \(\tilde{G}\)-invariant support \(\tau\)-tilting \(kG\)-module \(M\) such that \(\tilde{M}=_{\add}\induc_{G}^{\tilde{G}}M\) by \Cref{p-extension}.
  Hence, \(\tilde{X}\) is a direct summand of \(\induc_{G}^{\tilde{G}}M\).
  Therefore, \(\tilde{X}\) has a vertex contained in a Sylow \(p\)-subgroup of \(\tilde{G}\) properly.
\end{proof}

\section{Preliminaries for the block version of the main results}
We recall the definition of blocks of group algebras. Let \(G\) be a finite group. The group algebra \(kG\) has a unique decomposition
\begin{equation}\label{block dec 2021-11-10 11:00:13}
  kG=B_0\times \cdots \times B_l
\end{equation}
into the direct product of indecomposable algebras.
We call each indecomposable direct product component \(B_i\) a block of \(kG\) and the decomposition above the block decomposition. We remark that any block \(B_i\) is a two-sided ideal of \(kG\).

For any indecomposable \(kG\)-module \(U\), there exists a unique block \(B_i\) of \(kG\) such that \(U=B_iU\) and \(B_jU=0\) for all \(j\neq i\). Then we say that \(U\) lies in the block \(B_i\) or simply \(U\) is a \(B_i\)-module.
We denote by \(B_0(kG)\) the principal block of \(kG\), in which the trivial \(kG\)-module \(k_G\) lies.

Let \(G\) be a normal subgroup of a finite group \(\tilde{G}\), \(B\) a block of \(kG\) and \(\tilde{B}\) a block of \(k\tilde{G}\).
We say that \(\tilde{B}\) covers \(B\) (or that \(B\) is covered by \(\tilde{B}\)) if \(1_B 1_{\tilde{B}}\neq 0\), where \(1_B\) and \(1_{\tilde{B}}\) mean the respective identity element of \(B\) and \(\tilde{B}\).
% We remark that \(1_B\) and \(1_{\tilde{B}}\) are central idempotents of \(kG\) and \(k\tilde{G}\), respectively.

\begin{proposition}[{See \cite[Theorem 15.1, Lemma 15.3]{MR860771}}]\label{Remark:cover}
  With the notation above, the following are equivalent:
  \begin{enumerate}
    \item The block \(\tilde{B}\) covers \(B\).\label{cover item 0 2022-12-21 16:50:21}
    \item There exists a non-zero \(\tilde{B}\)-module \(\tilde{U}\) such that \(\restr_G^{\tilde{G}} \tilde{U}\) has a non-zero direct summand lying in \(B\).\label{Remark:cove direct summand as bimodule 2022-12-07 17:53:28}
    \item For any non-zero \(\tilde{B}\)-module \(\tilde{U}\), there exists a non-zero direct summand of \(\restr_G^{\tilde{G}} \tilde{U}\) lying in \(B\).\label{rem cov item 3 2022-12-21 15:52:39}
    \item For any non-zero \(\tilde{B}\)-module \(\tilde{U}\) and indecomposable direct summand \(V\) of \(\restr_{G}^{\tilde{G}}\tilde{U}\), there exists \(\tilde{g}\in \tilde{G}\) such that \(V\) lies in the block \(\tilde{g}B\tilde{g}^{-1}\). \label{ind res summand 2022-12-07 17:05:41 item 0}
    \item The block \(B\) is a direct summand of \(\tilde{B}\) as a \((kG,kG)\)-bimodule.\label{Remark:cove direct summand as bimodule 2022-12-07 16:36:27}
    \item The block \(\tilde{B}\) is a direct summand of \(k\tilde{G}B\tilde{G}\) as a \((k\tilde{G},k\tilde{G})\)-bimodule.
  \end{enumerate}
\end{proposition}

We denote by \(\inertiagp_{\tilde{G}}(B)\) the  inertial group of \(B\) in \(\tilde{G}\), that is
\begin{equation}
  \inertiagp_{\tilde{G}}(B):=\left\{ \tilde{g} \in \tilde{G} \setmid \tilde{g}B\tilde{g}^{-1} = B \right\}.
\end{equation}

\begin{remark}\label{only covered 2022-12-19 20:33:14}
  For a block \(\tilde{B}\) of \(k\tilde{G}\) and a block \(B\) of \(kG\), the block \(\tilde{B}\) covers only \(B\) if and only if \(\inertiagp_{\tilde{G}}(B)=\tilde{G}\) by \cite[Theorem 15.1 (1)]{MR860771}.
  Since \(\restr_{G}^{\tilde{G}}k_{\tilde{G}}\cong k_{G}\), the principal block \(B_0(kG)\) of \(kG\) is the only block of \(kG\) covered by the principal block \(B_0(k\tilde{G})\) of \(k\tilde{G}\) by the equivalence of \Cref{Remark:cover} \ref{cover item 0 2022-12-21 16:50:21}, \ref{rem cov item 3 2022-12-21 15:52:39}.
  Therefore, we have that \(\inertiagp_{\tilde{G}}(B_0(kG))=\tilde{G}\).
\end{remark}
\begin{proposition}\label{ind res summand 2022-12-07 17:05:41}
  Let \(G\) be a normal subgroup of a finite group \(\tilde{G}\), \(B\) a block of \(kG\) and \(U\) an indecomposable \(B\)-module.
  Then the following hold:
  \begin{enumerate}
    \item For a block \(\tilde{B}\) of \(k\tilde{G}\) covering the block \(B\), the module \(\restr_{G}^{\tilde{G}}\tilde{B}\induc_G^{\tilde{G}}U\) has a direct summand isomorphic to \(U\). In particular, the \(\tilde{B}\)-module \(\tilde{B}\induc_{G}^{\tilde{G}} U\) is non-zero. \label{ind res summand 2022-12-07 17:05:41 item1}
    \item  Any indecomposable direct summand \(\tilde{V}\) of \(\induc_G^{\tilde{G}}U\) lies in a block of \(k\tilde{G}\) covering \(B\).\label{ind res summand 2022-12-07 17:05:41 item2}
  \end{enumerate}
\end{proposition}
\begin{proof}
  Let \(U\) be an indecomposable \(B\)-module.
  By the equivalence of \Cref{Remark:cover} \ref{cover item 0 2022-12-21 16:50:21}, \ref{Remark:cove direct summand as bimodule 2022-12-07 16:36:27}, the block \(\tilde{B}\) has a direct summand \(B\) as a \((kG,kG)\)-bimodule.
  Hence, there exists a \((kG,kG)\)-bimodule \(B'\) such that \(\tilde{B}\cong B\oplus B'\) as a \((kG,kG)\)-bimodule.
  Therefore, we have that
  \begin{equation}
    \restr_{G}^{\tilde{G}}\tilde{B}\induc_G^{\tilde{G}}U
    \cong \restr_{G}^{\tilde{G}}\tilde{B}(k\tilde{G}\otimes_{kG}U)
    \cong \restr_{G}^{\tilde{G}}\tilde{B}\otimes_{kG} U
    \cong (B\oplus B')\otimes_{kG}U
    \cong U\oplus (B'\otimes_{kG} U),
  \end{equation}
  which prove \ref{ind res summand 2022-12-07 17:05:41 item1}.

  Let \(\tilde{V}\) be an indecomposable direct summand of \(\induc_{G}^{\tilde{G}}U\) lying in a block \(\tilde{A}\) of \(k\tilde{G}\).
  Since the restricted module \(\restr_{G}^{\tilde{G}}\tilde{V}\) is a direct summand of the \(kG\)-module \(\restr_{G}^{\tilde{G}}\induc_{G}^{\tilde{G}}U\), we have that the block \(\tilde{A}\) covers \(B\) by \Cref{mackey} and the equivalences of \Cref{Remark:cover} \ref{cover item 0 2022-12-21 16:50:21}, \ref{Remark:cove direct summand as bimodule 2022-12-07 17:53:28}, \ref{ind res summand 2022-12-07 17:05:41 item 0}.
  Hence, we get that \ref{ind res summand 2022-12-07 17:05:41 item2}.
\end{proof}
The following is a generalization of \cite[Corollary 5.5.6]{MR998775} (or \cite[Corollary 9.9.6]{MR1632299}).
\begin{proposition}\label{cover unique}
  Let \(G\) be a normal subgroup of a finite group \(\tilde{G}\) and \(B\) a block of \(kG\).
  If there exists an indecomposable \(B\)-module \(X\) such that \(\induc_G^{\tilde{G}}X\) is an indecomposable \(k\tilde{G}\)-module, then there exists only one block of \(k\tilde{G}\) covering \(B\).
\end{proposition}
\begin{proof}
  Let \(\tilde{A}\) and \(\tilde{B}\) be a block of \(k\tilde{G}\) covering \(B\).
  The modules \(\tilde{A}\induc_{G}^{\tilde{G}}X\)  and \(\tilde{B}\induc_{G}^{\tilde{G}}X\) are non-zero direct summands of the indecomposable \(k\tilde{G}\)-module \(\induc_{G}^{\tilde{G}}X\) by \Cref{ind res summand 2022-12-07 17:05:41} \ref{ind res summand 2022-12-07 17:05:41 item1}.
  Hence, we get that \(\tilde{B}\induc_{G}^{\tilde{G}}X\cong\tilde{A}\induc_{G}^{\tilde{G}}X\cong \induc_{G}^{\tilde{G}}X\) by the indecomposability of \(\induc_{G}^{\tilde{G}}X\).
  Since the non-zero \(k\tilde{G}\)-module \(\induc_{G}^{\tilde{G}}X\) lies in the blocks \(\tilde{A}\) and \(\tilde{B}\), we get that \(\tilde{A}=\tilde{B}\).
\end{proof}
\begin{corollary}[{See \cite[Corollary 5.5.6]{MR998775} or \cite[Corollary 9.9.6]{MR1632299}}]\label{covering block p-power 2022-12-12 17:59:37}
  If \(\tilde{G}/G\) is a \(p\)-group, then there exists only one block of \(k\tilde{G}\) covering \(B\).
\end{corollary}
\begin{proof}
  It immediately follows from \Cref{cover unique} and Green's indecomposability theorem (for example, see \cite{MR860771, MR131454, MR998775}).
\end{proof}

\begin{proposition}[{See \cite[Theorem 6.8.3]{MR3821517} or \cite[Theorem 5.5.10, Theorem 5.5.12]{MR998775}}]\label{Morita equivalence covering block}
  Let \(G\) be a normal subgroup of a finite group \(\tilde{G}\) and \(B\) a block of \(kG\).
  Then the following hold:
  \begin{enumerate}
    \item For any  block \(\beta\) of \(k\inertiagp_{\tilde{G}}(B)\) covering \(B\), there exists a block \(\tilde{B}\) of \(k\tilde{G}\) such that
          \begin{equation}
            \sum_{x\in [\tilde{G}/\inertiagp_{\tilde{G}}(B)]}x1_\beta x^{-1}=1_{\tilde{B}},
          \end{equation}
          and then \(\tilde{B}\) covers \(B\).
          Moreover, the correspondence sending \(\beta\) to \(\tilde{B}\) induces a bijection between the set of blocks of \(k\inertiagp_{\tilde{G}}(B)\) covering \(B\) and those of \(k\tilde{G}\) covering \(B\).\label{fong correspondence 2022-11-21 17:16:57}
    \item If \(\tilde{B}\) corresponds to \(\beta\) under the bijection of \ref{fong correspondence 2022-11-21 17:16:57}, then the induction functor
          \begin{equation}
            \begin{tikzcd}
              \induc_{\inertiagp_{\tilde{G}}(B)}^{\tilde{G}} \colon k\inertiagp_{\tilde{G}}(B)\lmod \ar[r] & k{\tilde{G}}\lmod
            \end{tikzcd}
          \end{equation}
          restricts to a Morita equivalence
          \begin{equation}\label{2020-03-25 15:17:09}
            \begin{tikzcd}
              \induc_{\inertiagp_{\tilde{G}}(B)}^{\tilde{G}} \colon \beta\lmod \ar[r]&  \tilde{B}\lmod
            \end{tikzcd}
          \end{equation}
          and its inverse functor is given by
          \begin{equation}\label{2022-11-21 17:18:47}
            \begin{tikzcd}
              \beta\restr_{\inertiagp_{\tilde{G}}(B)}^{\tilde{G}} \colon \tilde{B}\lmod \ar[r]&  \beta\lmod.
            \end{tikzcd}
          \end{equation}
          \label{Morita item 2022-06-09 01:35:47}
  \end{enumerate}
\end{proposition}
\begin{proposition}\label{fong conseq 2022-12-12 17:51:31}
  Let \(G\) be a normal subgroup of a finite group \(\tilde{G}\), \(B\) a block of \(kG\), \(U\) a \(B\)-module, \(\beta\) a block of \(k\inertiagp_{\tilde{G}}(B)\) covering \(B\) and \(\tilde{B}\) a block of \(k\tilde{G}\) covering \(B\) such that
  \begin{equation}
    \sum_{x\in [\tilde{G}/\inertiagp_{\tilde{G}}(B)]}x1_\beta x^{-1}=1_{\tilde{B}}.
  \end{equation}
  Then \(\tilde{B}\induc_G^{\tilde{G}}U\cong \induc_{\inertiagp_{\tilde{G}}(B)}^{\tilde{G}}\beta \induc_G^{\inertiagp_{\tilde{G}}(B)}U\).
\end{proposition}
\begin{proof}
  Let \(\tilde{B_1}=\tilde{B}, \ldots, \tilde{B_e}\) be the all blocks of \(k\tilde{G}\) covering \(B\).
  By \Cref{Morita equivalence covering block}, we can take \(\beta_1=\beta, \ldots, \beta_e\) the blocks of \(k\inertiagp_{\tilde{G}}(B)\) satisfying the induction functor \(\induc_{\inertiagp_{\tilde{G}}(B)}^{\tilde{G}}\) restricts to a Morita equivalence

  \begin{equation}
    \begin{tikzcd}
      \induc_{\inertiagp_{\tilde{G}}(B)}^{\tilde{G}} \colon \beta_i \lmod \ar[r]&\tilde{B_i}\lmod
    \end{tikzcd}
  \end{equation}
  for any \(i=1,\ldots, e\).
  By \Cref{ind res summand 2022-12-07 17:05:41} \ref{ind res summand 2022-12-07 17:05:41 item2}, we get the following isomorphism:
  \begin{equation}
    \induc_G^{\inertiagp_{\tilde{G}}(B)}U\cong \beta_1 \induc_G^{\inertiagp_{\tilde{G}}(B)}U\oplus \cdots \oplus \beta_e \induc_G^{\inertiagp_{\tilde{G}}(B)}U.
  \end{equation}
  Moreover, by \Cref{Theorem: Frobenius and projective} \ref{transitive induc 2022-06-09 02:03:23}, we have that
  \begin{equation}
    \induc_G^{\tilde{G}}U\cong\induc_{\inertiagp_{\tilde{G}}(B)}^{\tilde{G}} \induc_G^{\inertiagp_{\tilde{G}}(B)} U\cong \induc_{\inertiagp_{\tilde{G}}(B)}^{\tilde{G}}\beta_1 \induc_{G}^{\inertiagp_{\tilde{G}}(B)}U\oplus \cdots \oplus \induc_{\inertiagp_{\tilde{G}}(B)}^{\tilde{G}}\beta_e \induc_{G}^{\inertiagp_{\tilde{G}}(B)}U.
  \end{equation}
  Since the \(k\tilde{G}\)-module \(\induc_{\inertiagp_{\tilde{G}}(B)}^{\tilde{G}}\beta_i \induc_{G}^{\inertiagp_{\tilde{G}}(B)}U\) lies in the block \(\tilde{B}_i\) for any \(i=1,\ldots, e\), we get that
  \begin{equation}
    \tilde{B}_i\induc_G^{\tilde{G}}U\cong\induc_{\inertiagp_{\tilde{G}}(B)}^{\tilde{G}}\beta_i \induc_{G}^{\inertiagp_{\tilde{G}}(B)}U.
  \end{equation}
  Therefore, we complete the proof.
\end{proof}
\section{Block version of main results}\label{section block main}
In this \namecref{section block main}, we give the block versions of our \namecref{main-thm-1} stated in \Cref{sec main block}.
Let \(\Lambda\) be a finite dimensional algebra.
For \(\Lambda\)-modules \(M\) and \(N\), we write \(M \leq_{\add} N\) if \(\add M\subset \add N\).
This relation is clearly reflexive and transitive.
Moreover, if \(M \leq_{\add} N\) and \(N \leq_{\add} M\) then \(M=_{\add}N\) for any \(\Lambda\)-modules \(M\) and \(N\).
The following is the special case of the block version of \Cref{main-thm-1}.
\begin{theorem}
  \label{main-thm-1 block}
  Let \(G\) be a normal subgroup of a finite group \(\tilde{G}\), \(B\) a block of \(kG\) satisfying \(\inertiagp_{\tilde{G}}(B)=\tilde{G}\), \(\tilde{B}\) a block of \(k\tilde{G}\) covering \(B\) and \(\tilde{M}\) a support \(\tau\)-tilting \(\tilde{B}\)-module.
  If it holds that \(\tilde{B}\induc_{G}^{\tilde{G}}\restr_{G}^{\tilde{G}}\tilde{M} \in \add \tilde{M}\)
  and \(\tilde{M}\) is relatively \(G\)-projective,
  then we have that \(\restr_{G}^{\tilde{G}}\tilde{M}\) is a support \(\tau\)-tilting \(B\)-module.
  Moreover, if \((\tilde{M}, \tilde{P})\) is a support \(\tau\)-tilting pair for \(\tilde{B}\)
  corresponding to \(\tilde{M}\), then
  \((\restr_{G}^{\tilde{G}}\tilde{M}, \restr_{G}^{\tilde{G}}\tilde{P})\) is
  a support \(\tau\)-tilting pair for \(B\)
  corresponding to \(\restr_{G}^{\tilde{G}}\tilde{M}\).
\end{theorem}
\begin{proof}
  Let \((\tilde{M}, \tilde{P})\) be a support \(\tau\)-tilting pair for \(\tilde{B}\)
  corresponding to the support \(\tau\)-tilting \(\tilde{B}\)-module \(\tilde{M}\).
  Our assumption that \(\inertiagp_{\tilde{G}}(B)=\tilde{G}\) means the block \(B\) is the only block of \(kG\) covered by \(\tilde{B}\) by \Cref{only covered 2022-12-19 20:33:14}.
  Hence, we have that the restricted modules \(\restr_{G}^{\tilde{G}}\tilde{M}\) and \(\restr_{G}^{\tilde{G}}\tilde{P}\) are \(B\)-modules by \Cref{Remark:cover} \ref{ind res summand 2022-12-07 17:05:41 item 0}.

  First, we show that \((\restr_{G}^{\tilde{G}}\tilde{M}, \restr_{G}^{\tilde{G}}\tilde{P})\)
  is a \(\tau\)-rigid pair for \(B\).
  Since the \(\tilde{B}\)-module \(\tilde{M}\) is a support \(\tau\)-tilting \(\tilde{B}\)-module, it is a \(\tau\)-rigid \(\tilde{B}\)-module.
  Hence, we have that \(\restr_{G}^{\tilde{G}}\tilde{M}\) is a \(\tau\)-rigid \(B\)-module by \Cref{res-rigid}.
  On the other hand, by \Cref{Theorem: Frobenius and projective} \ref{frob adjoint 2022-12-26 19:17:16} we have that
  \begin{align}
    \Hom_{B}(\restr_{G}^{\tilde{G}}\tilde{P}, \restr_{G}^{\tilde{G}}\tilde{M})
     & \cong \Hom_{kG}(\restr_{G}^{\tilde{G}}\tilde{P}, \restr_{G}^{\tilde{G}}\tilde{M})                  \\
     & \cong \Hom_{k\tilde{G}}(\tilde{P}, \induc_{G}^{\tilde{G}}\restr_{G}^{\tilde{G}}\tilde{M})          \\
     & \cong \Hom_{\tilde{B}}(\tilde{P}, \tilde{B}\induc_{G}^{\tilde{G}}\restr_{G}^{\tilde{G}}\tilde{M}).
  \end{align}
  Now, by the assumption that \(\tilde{B}\induc_{G}^{\tilde{G}}\restr_{G}^{\tilde{G}}\tilde{M} \in \add \tilde{M}\),
  we have that \(\Hom_{\tilde{B}}(\tilde{P}, \tilde{B}\induc_{G}^{\tilde{G}}\restr_{G}^{\tilde{G}}\tilde{M})=0\)
  because \((\tilde{M}, \tilde{P})\) is a support \(\tau\)-tilting pair for \(\tilde{B}\).
  Therefore, we conclude that
  \((\restr_{G}^{\tilde{G}}\tilde{M}, \restr_{G}^{\tilde{G}}\tilde{P})\)
  is a \(\tau\)-rigid pair for \(B\).

  Next, we show that the \(\tau\)-rigid pair \((\restr_{G}^{\tilde{G}}\tilde{M}, \restr_{G}^{\tilde{G}}\tilde{P})\)
  is a support \(\tau\)-tilting pair for \(B\).
  We show that \(X\in \add\restr_G^{\tilde{G}}\tilde{M}\)
  under the assumption that
  \begin{equation}
    \Hom_{B}(X, \tau(\restr_G^{\tilde{G}}\tilde{M}))=\Hom_{B}(\restr_G^{\tilde{G}}\tilde{M}, \tau X)=\Hom_{B}(\restr_G^{\tilde{G}}\tilde{P}, X)=0,
  \end{equation}
  which implies that the pair \((\restr_{G}^{\tilde{G}}\tilde{M}, \restr_{G}^{\tilde{G}}\tilde{P})\) is a support \(\tau\)-tilting pair for \(B\) by \Cref{AIR Corollary2.13}.
  Under these assumptions, we have the following:
  \begin{align}
    \Hom_{\tilde{B}}(\tilde{B}\induc_{G}^{\tilde{G}}X, \tau\tilde{M})
     & \cong\Hom_{k\tilde{G}}(\induc_{G}^{\tilde{G}}X, \tau\tilde{M}) & \text{(\(\tau \tilde{M}\) is a \(\tilde{B}\)-module)}                                    \\
     & \cong\Hom_{kG}(X, \restr_{G}^{\tilde{G}}(\tau\tilde{M}))       & \text{(\Cref{Theorem: Frobenius and projective} \ref{frob adjoint 2022-06-02 15:44:10})} \\
     & \cong \Hom_{kG}(X, \tau(\restr_{G}^{\tilde{G}}\tilde{M}))      & \text{(\Cref{tau-restriction})}                                                          \\
     & \cong \Hom_{B}(X, \tau(\restr_{G}^{\tilde{G}}\tilde{M}))       & \text{(\(X\) and \(\tau(\restr_{G}^{\tilde{G}}\tilde{M})\) are \(B\)-modules)}           \\
     & =0.
  \end{align}
  \begin{align}
    \Hom_{\tilde{B}}(\tilde{M}, \tau(\tilde{B}\induc_{G}^{\tilde{G}}X))
     & \cong \Hom_{\tilde{B}}(\tilde{M}, \tilde{B}\induc_{G}^{\tilde{G}}(\tau X)) & \text{(\Cref{commute})}                                                                  \\
     & \cong \Hom_{k\tilde{G}}(\tilde{M}, \induc_{G}^{\tilde{G}}(\tau X))         & \text{(\(\tilde{M}\) is the \(\tilde{B}\)-module)}                                       \\
     & \cong \Hom_{kG}(\restr_{G}^{\tilde{G}}\tilde{M}, \tau X)                   & \text{(\Cref{Theorem: Frobenius and projective} \ref{frob adjoint 2022-12-26 19:17:16})} \\
     & \cong \Hom_{B}(\restr_{G}^{\tilde{G}}\tilde{M}, \tau X)                    & \text{(\(\restr_{G}^{\tilde{G}}\tilde{M}\) and \(\tau X\) are \(B\)-modules)}            \\
     & =0.
  \end{align}
  \begin{align}
    \Hom_{\tilde{B}}(\tilde{P}, \tilde{B}\induc_{G}^{\tilde{G}}X)
     & \cong \Hom_{k\tilde{G}}(\tilde{P}, \induc_{G}^{\tilde{G}}X) & \text{(\(\tilde{P}\) is a \(\tilde{B}\)-module)}                                         \\
     & \cong \Hom_{kG}(\restr_{G}^{\tilde{G}}\tilde{P}, X)         & \text{(\Cref{Theorem: Frobenius and projective} \ref{frob adjoint 2022-12-26 19:17:16})} \\
     & \cong \Hom_{B}(\restr_{G}^{\tilde{G}}\tilde{P}, X)          & \text{(\(\restr_{G}^{\tilde{G}}\tilde{P}\) and \(X\) are \(B\)-modules)}                 \\
     & =0.
  \end{align}

  By these three isomorphisms and the fact that \((\tilde{M}, \tilde{P})\)
  is a support \(\tau\)-tilting pair for \(\tilde{B}\),
  applying \Cref{AIR Corollary2.13},
  we have that \(\tilde{B}\induc_{G}^{\tilde{G}}X\in \add \tilde{M}\).
  Also, since the block \(\tilde{B}\) covers \(B\), the \(B\)-module \(X\) is a direct summand of \( \restr_{G}^{\tilde{G}}\tilde{B}\induc_{G}^{\tilde{G}}X\) by \Cref{ind res summand 2022-12-07 17:05:41} \ref{ind res summand 2022-12-07 17:05:41 item1}.
  Therefore, we have that \(X\in \add \restr_{G}^{\tilde{G}}\tilde{M}\).
\end{proof}

We consider equivalent conditions to the assumption of \Cref{main-thm-1 block}.
First, we give the \namecrefs{dec lemma block 2022-12-07 11:51:53} which can be applied in case of rigid \(\tilde{B}\)-modules
not only support \(\tau\)-tilting \(\tilde{B}\)-modules.
The following \namecref{dec lemma block 2022-12-07 11:51:53} is the block version of \Cref{dec lemma group algebra 2022-12-07 11:51:53}, which is helpful to prove \Cref{main-thm-2 block}.

\begin{lemma}\label{dec lemma block 2022-12-07 11:51:53}
  Let \(\tilde{G}\) be a finite group, \(\tilde{B}\) a block of \(k\tilde{G}\), \(\tilde{M}\) a rigid \(\tilde{B}\)-module and \(L\) a \(k\tilde{G}\)-module.
  If it holds that \(\tilde{B}(S\otimes_k \tilde{M})\in \add \tilde{M}\) for any composition factor \(S\) of \(L\),
  then the following isomorphism as \(\tilde{B}\)-modules holds:
  \begin{equation}
    \tilde{B}(L\otimes_k \tilde{M}) \cong \bigoplus_{S}\tilde{B}(S \otimes_k \tilde{M}),
  \end{equation}
  where \(S\) runs over all composition factors of \(L\).
\end{lemma}
\begin{proof}
  A similar proof of \Cref{dec lemma group algebra 2022-12-07 11:51:53} works in this setting.
  Let \(L\) be an arbitrarily \(k\tilde{G}\)-module and \(\tilde{M}\) a rigid \(\tilde{B}\)-module satisfying that
  \begin{equation}\label{dec lemma block 2022-12-07 11:51:53 assumption}
    \tilde{B}(S\otimes_k \tilde{M})\in \add \tilde{M} \text{ for any composition factors \(S\) of \(L\).}
  \end{equation}
  We use induction on the composition length \(\ell(L)\) of \(L\).
  If \(\ell(L)=1\), there is nothing to prove.
  Hence, we assume that \(\ell(L)\geq 2\) and that the statement for any \(k\tilde{G}\)-module \(L'\) satisfying \(\ell(L')<\ell(L)\) is true.
  Let \(T\) be a simple submodule of \(L\).
  We get the exact sequence
  \begin{equation}\label{ex_seq tensor 2022-12-07 13:58:32}
    \begin{tikzcd}
      0\ar[r]&\tilde{B}(T\otimes_k \tilde{M})\ar[r]&\tilde{B}(L\otimes_k \tilde{M})\ar[r]&\tilde{B}((L/T)\otimes_k \tilde{M})\ar[r]&0
    \end{tikzcd}
  \end{equation}
  obtained by applying the exact functor \(\tilde{B}(-\otimes_k \tilde{M})\)
  to the exact sequence
  \begin{equation}
    \begin{tikzcd}
      0\ar[r]&T\ar[r]&L\ar[r]&L/T\ar[r]&0.
    \end{tikzcd}
  \end{equation}
  By the rigidity of \(\tilde{M}\), the assumption \eqref{dec lemma block 2022-12-07 11:51:53 assumption} and the assumption of the induction, the sequence \eqref{ex_seq tensor 2022-12-07 13:58:32} splits, and we get that
  \begin{equation}
    \tilde{B}(L\otimes_k \tilde{M})\cong \tilde{B}(T\otimes_k \tilde{M})\oplus \tilde{B}((L/T)\otimes_k \tilde{M}) \cong \tilde{B}(T\otimes_k \tilde{M})\oplus \bigoplus_{S'}\tilde{B}(S'\otimes_k \tilde{M})\cong \bigoplus_{S}\tilde{B}(S\otimes_k \tilde{M}),
  \end{equation}
  where \(S'\) and \(S\) run over all composition factors of \(L/T\) and \(L\), respectively.
\end{proof}

\begin{lemma}\label{equivalence-condition block}
  Let \(G\) be a normal subgroup of a finite group \(\tilde{G}\), \(\tilde{B}\) a block of \(k\tilde{G}\) and \(\tilde{M}\) a rigid \(\tilde{B}\)-module.
  Then the following conditions are equivalent:
  \begin{enumerate}
    \item \(\tilde{B}\induc_{G}^{\tilde{G}}\restr_{G}^{\tilde{G}} \tilde{M} \in \add \tilde{M}\) and \(\tilde{M}\) is relatively \(G\)-projective.\label{equivalence-condition block item1}
    \item \(\tilde{B}(S\otimes_k \tilde{M}) \in \add \tilde{M}\) for each simple \(k(\tilde{G}/G)\)-module \(S\).\label{equivalence-condition block item2}
  \end{enumerate}
\end{lemma}
\begin{proof}
  By \Cref{Theorem: Frobenius and projective}, we have that
  \begin{equation}
    \tilde{B}\induc_{G}^{\tilde{G}}\restr_{G}^{\tilde{G}} \tilde{M}
    \cong \tilde{B}\induc_{G}^{\tilde{G}} (k_{G} \otimes_k \restr_{G}^{\tilde{G}}\tilde{M})
    \cong \tilde{B}((\induc_{G}^{\tilde{G}} k_{G}) \otimes_k \tilde{M})
    \cong \tilde{B}(k(\tilde{G}/G) \otimes_k \tilde{M}).
  \end{equation}
  \ref{equivalence-condition block item1} \(\Rightarrow\) \ref{equivalence-condition block item2}.
  By the assumptions, we have that \(\tilde{B}\induc_{G}^{\tilde{G}}\restr_{G}^{\tilde{G}} \tilde{M}=_{\add}\tilde{M}\).
  Hence, by \Cref{Theorem: Frobenius and projective} we get that
  \begin{align}
    \tilde{B}(S\otimes_k \tilde{M})
     & =_{\add} \tilde{B}(S\otimes_k \tilde{B}\induc_{G}^{\tilde{G}}\restr_{G}^{\tilde{G}}\tilde{M})            \\
     & \leq_{\add} \tilde{B}(S\otimes_k \induc_{G}^{\tilde{G}}\restr_{G}^{\tilde{G}}\tilde{M})                  \\
     & \cong \tilde{B}\induc_{G}^{\tilde{G}}(\restr_{G}^{\tilde{G}}S\otimes_k\restr_{G}^{\tilde{G}}\tilde{M})   \\
     & \cong \tilde{B}\induc_{G}^{\tilde{G}}( k_{G}^{\oplus \dim_k S} \otimes_k\restr_{G}^{\tilde{G}}\tilde{M}) \\
     & =_{\add} \tilde{B}\induc_{G}^{\tilde{G}}\restr_{G}^{\tilde{G}} \tilde{M}                                 \\
     & =_{\add} \tilde{M},
  \end{align}
  for any simple \(k(\tilde{G}/G)\)-module \(S\),
  which implies that  \(\tilde{B}(S\otimes_k \tilde{M})\in \add \tilde{M}\).

  \noindent
  \ref{equivalence-condition block item2} \(\Rightarrow\) \ref{equivalence-condition block item1}.
  By \Cref{dec lemma block 2022-12-07 11:51:53}, we have that
  \begin{equation}
    \tilde{B}\induc_{G}^{\tilde{G}}\restr_{G}^{\tilde{G}} \tilde{M}
    \cong \tilde{B}(k(\tilde{G}/G) \otimes_k \tilde{M})
    \cong \bigoplus_{S}\tilde{B}(S \otimes_k \tilde{M}),
  \end{equation}
  where \(S\) runs over all composition factors of the \(k\tilde{G}\)-module \(k(\tilde{G}/G)\).
  Therefore, the assumption implies that \(\tilde{B}\induc_{G}^{\tilde{G}}\restr_{G}^{\tilde{G}} \tilde{M}\in \add \tilde{M}\).
  Moreover, since the trivial \(k\tilde{G}\)-module \(k_{\tilde{G}}\) appears as a composition factor of \(k(\tilde{G}/G)\), we have that the \(\tilde{B}\)-module \(\tilde{M}\) appears as a direct summand of \(\tilde{B}\induc_{G}^{\tilde{G}}\restr_{G}^{\tilde{G}} \tilde{M}\), that is \(\tilde{M}\) is a relatively \(G\)-projective \(k\tilde{G}\)-module.
\end{proof}

We give the equivalent conditions to the assumption of \Cref{main-thm-1 block}.

\begin{theorem}\label{main-thm-2 block}
  Let \(G\) be a normal subgroup of a finite group, \(B\) a block of \(kG\) satisfying \(\inertiagp_{\tilde{G}}(B)=\tilde{G}\) and \(\tilde{B}\) a block of \(k\tilde{G}\) covering \(B\).
  Let \(\tilde{M}\) be a support \(\tau\)-tilting \(\tilde{B}\)-module.
  Then the following conditions are equivalent:
  \begin{enumerate}
    \item \(\tilde{M}=_{\add} \tilde{B}\induc_{G}^{\tilde{G}} M\) for some \(\tilde{G}\)-invariant support \(\tau\)-tilting \(B\)-module \(M\).\label{main thm block item1 2022-12-13 14:43:45}
    \item \(\tilde{B}\induc_{G}^{\tilde{G}}\restr_{G}^{\tilde{G}} \tilde{M} \in \add \tilde{M}\) and \(\tilde{M}\) is relatively \(G\)-projective.\label{main thm block item2 2022-12-13 14:43:45}
    \item \(\tilde{B}(S\otimes_k \tilde{M}) \in \add \tilde{M}\) for each simple \(k(\tilde{G}/G)\)-module \(S\).\label{main thm block item3 2022-12-13 14:43:45}
  \end{enumerate}
\end{theorem}
\begin{proof}
  \ref{main thm block item1 2022-12-13 14:43:45} \(\Rightarrow\) \ref{main thm block item2 2022-12-13 14:43:45}.
  Assume that \(\tilde{M}=_{\add} \tilde{B}\induc_{G}^{\tilde{G}} M\) for some \(\tilde{G}\)-invariant
  support \(\tau\)-tilting \(B\)-module \(M\).
  Then clearly \(\tilde{M}\) is a relatively \(G\)-projective \(\tilde{B}\)-module, and we get that
  \begin{align}
    \tilde{B}\induc_{G}^{\tilde{G}}\restr_{G}^{\tilde{G}} \tilde{M}
     & =_{\add}\tilde{B}\induc_{G}^{\tilde{G}}\restr_{G}^{\tilde{G}}\tilde{B} \induc_{G}^{\tilde{G}} M \\
     & \leq_{\add}\tilde{B}\induc_{G}^{\tilde{G}}\restr_{G}^{\tilde{G}}\induc_{G}^{\tilde{G}} M        \\
     & \cong \tilde{B}\induc_{G}^{\tilde{G}}(\bigoplus_{\tilde{g}\in [\tilde{G}/G]}\tilde{g}M)         \\
     & \cong \bigoplus_{\tilde{g}\in [\tilde{G}/G]}\tilde{B}\induc_{G}^{\tilde{G}} M                   \\
     & =_{\add}\tilde{M}.
  \end{align}
  Hence, we get \(\tilde{B}\induc_{G}^{\tilde{G}}\restr_{G}^{\tilde{G}} \tilde{M}\in \add \tilde{M}\).

  \noindent
  \ref{main thm block item2 2022-12-13 14:43:45} \(\Rightarrow\) \ref{main thm block item1 2022-12-13 14:43:45}.
  Assume that \(\tilde{B}\induc_{G}^{\tilde{G}}\restr_{G}^{\tilde{G}}\tilde{M} \in \add \tilde{M}\) and that \(\tilde{M}\) is relatively \(G\)-projective.
  Put \(M:=\restr_{G}^{\tilde{G}}\tilde{M}\).
  Then by \Cref{restriction-invariance}, \Cref{Remark:cover} \ref{ind res summand 2022-12-07 17:05:41 item 0}, \Cref{only covered 2022-12-19 20:33:14} and \Cref{main-thm-1 block}, \(M\) is a \(\tilde{G}\)-invariant support \(\tau\)-tilting \(B\)-module.
  We show that \(\tilde{B}\induc_{G}^{\tilde{G}} M =_{\add}\tilde{M}\), that is \(\add(\tilde{B}\induc_{G}^{\tilde{G}} M) =\add\tilde{M}\).
  By the assumption that \(\tilde{B}\induc_{G}^{\tilde{G}}M =\tilde{B}\induc_{G}^{\tilde{G}}\restr_{G}^{\tilde{G}}\tilde{M}\in \add \tilde{M}\),
  we have that \(\add(\tilde{B}\induc_{G}^{\tilde{G}} M) \subset\add\tilde{M}\).
  On the other hand, since \(\tilde{M}\) is relatively \(G\)-projective,
  \(\tilde{M}\) is a direct summand of \(\induc_{G}^{\tilde{G}}\restr_{G}^{\tilde{G}} \tilde{M}=\induc_{G}^{\tilde{G}} M\).
  Moreover, since \(\tilde{M}\) lies in \(\tilde{B}\), \(\tilde{M}\) is a direct summand of \(\tilde{B}\induc_{G}^{\tilde{G}} M\).
  Hence, we have \(\add\tilde{M} \subset \add(\tilde{B}\induc_{G}^{\tilde{G}} M)\).

  \noindent
  \ref{main thm block item2 2022-12-13 14:43:45} \(\Leftrightarrow\) \ref{main thm block item3 2022-12-13 14:43:45}
  Since support \(\tau\)-tilting \(\tilde{B}\)-modules
  are rigid \(\tilde{B}\)-modules, the equivalence follows from \Cref{equivalence-condition block}.
\end{proof}

\begin{corollary}\label{cor1 block}
  Let \(G\) be a normal subgroup of a finite group \(\tilde{G}\), \(B\) a block of \(kG\) satisfying \(\inertiagp_{\tilde{G}}(B)=\tilde{G}\) and \(\tilde{B}\) a block of \(k\tilde{G}\) covering \(B\).
  We denote by \((\stautilt B)^{\tilde{G}}\) the subset of \(\stautilt B\) consisting of \(\tilde{G}\)-invariant support \(\tau\)-tilting \(B\)-modules and by \((\stautilt \tilde{B})^{\star \star}\) the subset of \(\stautilt \tilde{B}\) consisting of support \(\tau\)-tilting \(\tilde{B}\)-modules satisfying the equivalent conditions of \Cref{main-thm-2 block}.
  Then the functor \(\tilde{B}\induc_{G}^{\tilde{G}}\) induces
  a poset isomorphism
  \begin{equation}\label{block well def iso 2022-12-13 16:30:18}
    \begin{tikzcd}[row sep=0pt]
      (\stautilt B)^{\tilde{G}}\ar[r,"\sim"] &(\stautilt \tilde{B})^{\star \star}\\
      M\ar[r,mapsto]&\tilde{B}\induc_{G}^{\tilde{G}}M.
    \end{tikzcd}
  \end{equation}
  In particular, the functor \(\tilde{B}\induc_{G}^{\tilde{G}}\) induces the poset monomorphism
  \begin{equation}\label{block well def map 2022-12-13 16:30:18}
    \begin{tikzcd}[row sep=0pt]
      (\stautilt B)^{\tilde{G}}\ar[r]& \stautilt \tilde{B}\\
      M\ar[r,mapsto]&\tilde{B}\induc_{G}^{\tilde{G}}M.
    \end{tikzcd}
  \end{equation}
\end{corollary}

\begin{proof}
  By \cite[Theorem 3.3]{https://doi.org/10.48550/arxiv.2208.14680},
  the map \eqref{block well def map 2022-12-13 16:30:18} is well-defined.
  Moreover, by the exactness of the functor \(\tilde{B}\induc_{G}^{\tilde{G}}\), if \(N\leq M\) in \(\stautilt B\) then \(\tilde{B}\induc_{G}^{\tilde{G}} N \leq \tilde{B}\induc_{G}^{\tilde{G}} M\) in \(\stautilt \tilde{B}\).
  Therefore, the map \eqref{block well def map 2022-12-13 16:30:18} is a poset homomorphism.
  % Also, the map is a poset homomorphism,
  % because for an epimorphism
  % \(
  % \begin{tikzcd}[column sep=20pt]
  %   M^{\oplus r}\ar[r,->>]& N
  % \end{tikzcd}
  % \),
  % the induced homomorphism
  % \(
  % \begin{tikzcd}[column sep=20pt]
  %   (\tilde{B}\induc_{G}^{\tilde{G}} M)^{\oplus r} \ar[r,->>]& \tilde{B}\induc_{G}^{\tilde{G}} N
  % \end{tikzcd}
  % \)
  % is an epimorphism, which means that if \(N\leq M\) then \(\tilde{B}\induc_{G}^{\tilde{G}} N \leq \tilde{B}\induc_{G}^{\tilde{G}} M\).

  It remains to show that the map \eqref{block well def map 2022-12-13 16:30:18} restricts to a poset isomorphism \eqref{block well def iso 2022-12-13 16:30:18}.
  By the definition of \((\stautilt \tilde{B})^{\star \star}\) and the above argument, the map \eqref{block well def iso 2022-12-13 16:30:18} is well-defined and a poset homomorphism.
  Also, for any relatively \(G\)-projective support \(\tau\)-tilting \(\tilde{B}\)-module \(\tilde{M}\) with \(\tilde{B}\induc_{G}^{\tilde{G}}\restr_{G}^{\tilde{G}}\tilde{M}\in \tilde{M}\),
  by \Cref{main-thm-2 block}, we can take a \(\tilde{G}\)-invariant support \(\tau\)-tilting \(B\)-module \(M\) satisfying \(\tilde{B}\induc_{G}^{\tilde{G}} M=_{\add} \tilde{M}\).
  Hence, the map \eqref{block well def iso 2022-12-13 16:30:18} is surjective.
  Also, assume that two \(\tilde{G}\)-invariant support \(\tau\)-tilting \(B\)-modules
  \(M\) and \(N\) satisfy that \(\tilde{B}\induc_{G}^{\tilde{G}} M =_{\add} \tilde{B}\induc_{G}^{\tilde{G}} N\).
  We have that
  \begin{equation}
    M
    \leq_{\add} \restr_{G}^{\tilde{G}}\tilde{B}\induc_{G}^{\tilde{G}}M
    \leq_{\add} \restr_{G}^{\tilde{G}}\induc_{G}^{\tilde{G}}M
    \cong \bigoplus_{\tilde{g}\in [\tilde{G}/G]} \tilde{g}M
    =_{\add} M
  \end{equation}
  by \Cref{ind res summand 2022-12-07 17:05:41} \ref{ind res summand 2022-12-07 17:05:41 item1}, \Cref{mackey} and the \(\tilde{G}\)-invariance of \(M\).
  Therefore, we have that \(M=_{\add} \restr_{G}^{\tilde{G}}\tilde{B}\induc_{G}^{\tilde{G}}M\) and that \(N=_{\add} \restr_{G}^{\tilde{G}}\tilde{B}\induc_{G}^{\tilde{G}}N\) similarly.
  Hence, we get that \(M=_{\add}\restr_{G}^{\tilde{G}}\tilde{B}\induc_{G}^{\tilde{G}}M=_{\add}\restr_{G}^{\tilde{G}}\tilde{B}\induc_{G}^{\tilde{G}}N=_{\add}N\), which implies that the map \eqref{block well def iso 2022-12-13 16:30:18} is injective.
  This completes the proof of the first assertion.

  The latter assertion immediately follows from the fact that the map \eqref{block well def map 2022-12-13 16:30:18} is the composition of the poset isomorphism \eqref{block well def iso 2022-12-13 16:30:18} and
  the inclusion map
  \(
  \begin{tikzcd}
    (\stautilt \tilde{B})^{\star \star}\ar[r,hook] &\stautilt \tilde{B}.
  \end{tikzcd}
  \)
\end{proof}
The following is the block version of \Cref{main-thm-1}.
\begin{theorem}
  \label{main-thm-1 block general}
  Let \(G\) be a normal subgroup of a finite group \(\tilde{G}\), \(B\) a block of \(kG\), \(\tilde{B}\) a block of \(k\tilde{G}\) covering \(B\), \(\beta\) the block of \(k\inertiagp_{\tilde{G}}(B)\) satisfying
  \begin{equation}
    \sum_{x\in [\tilde{G}/\inertiagp_{\tilde{G}}(B)]}x1_\beta x^{-1}=1_{\tilde{B}}
  \end{equation}
  and \(\tilde{M}\) a support \(\tau\)-tilting \(\tilde{B}\)-module.
  If it holds that \(\beta\induc_{G}^{\inertiagp_{\tilde{G}}(B)}\restr_{G}^{\inertiagp_{\tilde{G}}(B)}\beta\restr_{\inertiagp_{\tilde{G}}(B)}^{\tilde{G}}\tilde{M} \in \add \beta\restr_{\inertiagp_{\tilde{G}}(B)}^{\tilde{G}}\tilde{M}\)
  and \(\beta\restr_{\inertiagp_{\tilde{G}}(B)}^{\tilde{G}}\tilde{M}\) is relatively \(G\)-projective,
  then we have that \(\restr_{G}^{\inertiagp_{\tilde{G}}(B)}\beta\restr_{\inertiagp_{\tilde{G}}(B)}^{\tilde{G}}\tilde{M}\) is a support \(\tau\)-tilting \(B\)-module.
  Moreover, if \((\tilde{M}, \tilde{P})\) is a support \(\tau\)-tilting pair for \(\tilde{B}\)
  corresponding to \(\tilde{M}\), then the pair
  \begin{equation}
    (\restr_{G}^{\inertiagp_{\tilde{G}}(B)}\beta\restr_{\inertiagp_{\tilde{G}}(B)}^{\tilde{G}}\tilde{M}, \restr_{G}^{\inertiagp_{\tilde{G}}(B)}\beta\restr_{\inertiagp_{\tilde{G}}(B)}^{\tilde{G}}\tilde{P})
  \end{equation}
  is a support \(\tau\)-tilting pair for \(B\)
  corresponding to \(\restr_{G}^{\inertiagp_{\tilde{G}}(B)}\beta\restr_{\inertiagp_{\tilde{G}}(B)}^{\tilde{G}}\tilde{M}\).
\end{theorem}
\begin{proof}
  Since the functor
  \begin{equation}\label{2022-12-12 16:27:08}
    \begin{tikzcd}
      \beta\restr_{\inertiagp_{\tilde{G}}(B)}^{\tilde{G}} \colon \tilde{B}\lmod \ar[r]&  \beta\lmod
    \end{tikzcd}
  \end{equation}
  is a Morita equivalence by \Cref{Morita equivalence covering block}, the module \(\beta\restr_{\inertiagp_{\tilde{G}}(B)}^{\tilde{G}}\tilde{M}\) is a support \(\tau\)-tilting \(\beta\)-module and
  \begin{equation}
    (\beta\restr_{\inertiagp_{\tilde{G}}(B)}^{\tilde{G}}\tilde{M}, \beta\restr_{\inertiagp_{\tilde{G}}(B)}^{\tilde{G}}\tilde{P})
  \end{equation}
  is a corresponding support \(\tau\)-tilting pair for \(\beta\).
  Hence, by \Cref{main-thm-1 block} it immediately follows the consequence.
\end{proof}
The following is the block version of \Cref{main-thm-2}.
\begin{theorem}\label{main-thm-2 block general}
  Let \(G\) be a normal subgroup of a finite group \(\tilde{G}\), \(B\) a block of \(kG\), \(\tilde{B}\) a block of \(k\tilde{G}\) covering \(B\), \(\beta\) the block of \(k\inertiagp_{\tilde{G}}(B)\) satisfying
  \begin{equation}
    \sum_{x\in [\tilde{G}/\inertiagp_{\tilde{G}}(B)]}x1_\beta x^{-1}=1_{\tilde{B}}
  \end{equation}
  and \(\tilde{M}\) a support \(\tau\)-tilting \(\tilde{B}\)-module.
  Then the following conditions are equivalent:
  \begin{enumerate}
    \item \(\tilde{M}=_{\add} \tilde{B}\induc_{G}^{\tilde{G}} M\) for some \(\inertiagp_{\tilde{G}}(B)\)-invariant support \(\tau\)-tilting \(B\)-module \(M\).\label{main-thm-2 block general item1}
    \item \(\beta\induc_{G}^{\inertiagp_{\tilde{G}}(B)}\restr_{G}^{\inertiagp_{\tilde{G}}(B)}\beta\restr_{\inertiagp_{\tilde{G}}(B)}^{\tilde{G}}\tilde{M} \in \add \beta\restr_{\inertiagp_{\tilde{G}}(B)}^{\tilde{G}}\tilde{M}\) and \(\beta\restr_{\inertiagp_{\tilde{G}}(B)}^{\tilde{G}}\tilde{M}\) is relatively \(G\)-projective.\label{main-thm-2 block general item2}
    \item \(\beta(S\otimes_k \beta\restr_{\inertiagp_{\tilde{G}}(B)}^{\tilde{G}}\tilde{M}) \in \add \beta\restr_{\inertiagp_{\tilde{G}}(B)}^{\tilde{G}}\tilde{M}\) for each simple \(k(\inertiagp_{\tilde{G}}(B)/G)\)-module \(S\).\label{main-thm-2 block general item3}
  \end{enumerate}
\end{theorem}
\begin{proof}
  We remark that the module \(\beta\restr_{\inertiagp_{\tilde{G}}(B)}^{\tilde{G}}\tilde{M}\) is a support \(\tau\)-tilting \(\beta\)-module since the functor
  \begin{equation}\label{2022-12-12 17:23:54}
    \begin{tikzcd}
      \beta\restr_{\inertiagp_{\tilde{G}}(B)}^{\tilde{G}} \colon \tilde{B}\lmod \ar[r]&  \beta\lmod
    \end{tikzcd}
  \end{equation}
  is a Morita equivalence by \Cref{Morita equivalence covering block}.

  \noindent
  \ref{main-thm-2 block general item1} \(\Rightarrow\) \ref{main-thm-2 block general item2}.
  Assume that \(\tilde{M}=_{\add} \tilde{B}\induc_{G}^{\tilde{G}} M\) for some \(\inertiagp_{\tilde{G}}(B)\)-invariant
  support \(\tau\)-tilting \(B\)-module \(M\).
  By \cite[Theorem 3.3]{https://doi.org/10.48550/arxiv.2208.14680}, the module \(\beta\induc_{G}^{\inertiagp_{\tilde{G}}(B)}M\) is a support \(\tau\)-tilting \(\beta\)-module.
  Since the functor
  \begin{equation}\label{2022-12-15 15:50:00}
    \begin{tikzcd}
      \induc_{\inertiagp_{\tilde{G}}(B)}^{\tilde{G}} \colon \beta\lmod \ar[r]&  \tilde{B}\lmod
    \end{tikzcd}
  \end{equation}
  is a Morita equivalence with the inverse functor \eqref{2022-12-12 17:23:54}.
  we have \(\induc_{\inertiagp_{\tilde{G}}(B)}^{\tilde{G}}\beta \restr_{\inertiagp_{\tilde{G}}(B)}^{\tilde{G}} \tilde{M} \cong \tilde{M}\).
  Also, by the assumption and \Cref{fong conseq 2022-12-12 17:51:31}, we get that \(\tilde{M}=_{\add}\tilde{B}\induc_{G}^{\tilde{G}}M\cong \induc_{\inertiagp_{\tilde{G}}(B)}^{\tilde{G}}\beta\induc_{G}^{\inertiagp_{\tilde{G}}(B)}M\).
  Therefore, we have that \(\induc_{\inertiagp_{\tilde{G}}(B)}^{\tilde{G}}\beta \restr_{\inertiagp_{\tilde{G}}(B)}^{\tilde{G}} \tilde{M}=_{\add} \induc_{\inertiagp_{\tilde{G}}(B)}^{\tilde{G}}\beta\induc_{G}^{\inertiagp_{\tilde{G}}(B)}M\).
  Hence, by the fact that the functor \eqref{2022-12-15 15:50:00} is a Morita equivalence again, we have that \(\beta \restr_{\inertiagp_{\tilde{G}}(B)}^{\tilde{G}} \tilde{M}=_{\add} \beta\induc_{G}^{\inertiagp_{\tilde{G}}(B)}M\).
  Therefore, we get the consequence \ref{main-thm-2 block general item2} by the equivalence of \Cref{main-thm-2 block}. \ref{main thm block item1 2022-12-13 14:43:45} and \ref{main thm block item2 2022-12-13 14:43:45}.

  \noindent
  \ref{main-thm-2 block general item2} \(\Rightarrow\) \ref{main-thm-2 block general item1}.
  Since \(\beta\restr_{\inertiagp_{\tilde{G}}(B)}^{\tilde{G}}\tilde{M}\) is a support \(\tau\)-tilting \(\beta\)-module, there exists an \(\inertiagp_{\tilde{G}}(B)\)-invariant support \(\tau\)-tilting \(B\)-module \(M\) such that \(\beta\restr_{\inertiagp_{\tilde{G}}(B)}^{\tilde{G}}\tilde{M}=_{\add} \beta\induc_{G}^{\inertiagp_{\tilde{G}}(B)}M\) by the assumptions and \Cref{main-thm-2 block}.
  Therefore, by \Cref{fong conseq 2022-12-12 17:51:31}, we get that
  \begin{equation}
    \tilde{M} \cong  \induc_{\inertiagp_{\tilde{G}}(B)}^{\tilde{G}}\beta\restr_{\inertiagp_{\tilde{G}}(B)}^{\tilde{G}}\tilde{M}  =_{\add} \induc_{\inertiagp_{\tilde{G}}(B)}^{\tilde{G}}\beta\induc_{G}^{\inertiagp_{\tilde{G}}(B)}M \cong \tilde{B}\induc_{G}^{\tilde{G}} M.
  \end{equation}

  \noindent
  \ref{main-thm-2 block general item2} \(\Leftrightarrow\) \ref{main-thm-2 block general item3}.
  Since the support \(\tau\)-tilting \(\beta\)-module \(\beta\restr_{\inertiagp_{\tilde{G}}(B)}^{\tilde{G}}\tilde{M}\)
  is the rigid \(\beta\)-module, the equivalence
  follows from \Cref{equivalence-condition block}.
\end{proof}
\begin{corollary}\label{cor1 block general}
  Let \((\stautilt B)^{\inertiagp_{\tilde{G}}(B)}\) be the subset of
  \(\stautilt B\) consisting of \(\inertiagp_{\tilde{G}}(B)\)-invariant support \(\tau\)-tilting
  \(B\)-modules
  and
  \((\stautilt \tilde{B})^{\star \star \star}\) the subset of
  \(\stautilt \tilde{B}\) consisting of support \(\tau\)-tilting
  \(\tilde{B}\)-modules satisfying the equivalent conditions of \Cref{main-thm-2 block general}.
  Then the functor \(\tilde{B}\induc_{G}^{\tilde{G}}\) induces
  a poset isomorphism
  \begin{equation}\label{iso block general 2022-12-15 18:23:37}
    \begin{tikzcd}[row sep=0pt]
      (\stautilt B)^{\inertiagp_{\tilde{G}}(B)}\ar[r]& (\stautilt \tilde{B})^{\star \star \star}\\
      M\ar[r,mapsto]&\tilde{B}\induc_{G}^{\tilde{G}} M.
    \end{tikzcd}
  \end{equation}
  In particular, the functor \(\tilde{B}\induc_{G}^{\tilde{G}}\) induces the poset monomorphism
  \begin{equation}\label{inj block general 2022-12-15 18:23:37}
    \begin{tikzcd}[row sep=0pt]
      (\stautilt B)^{\inertiagp_{\tilde{G}}(B)}\ar[r]& \stautilt \tilde{B}\\
      M\ar[r,mapsto]&\tilde{B}\induc_{G}^{\tilde{G}} M.
    \end{tikzcd}
  \end{equation}
\end{corollary}
\begin{proof}
  Let \((\stautilt \beta)^{\star \star}\) be the subset of \(\stautilt \beta\) consisting of support \(\tau\)-tilting \(\beta\)-modules satisfying the equivalent conditions of \Cref{main-thm-2 block}.
  Since the functor
  \begin{equation}\label{2022-12-12 17:40:47}
    \begin{tikzcd}
      \induc_{\inertiagp_{\tilde{G}}(B)}^{\tilde{G}} \colon \beta\lmod \ar[r]&  \tilde{B}\lmod
    \end{tikzcd}
  \end{equation}
  is a Morita equivalence, we have poset isomorphisms
  \begin{equation}
    \begin{tikzcd}[row sep=0pt]
      \stautilt \beta \ar[r]&\stautilt \tilde{B}\\
      M\ar[r,mapsto]&\induc_{\inertiagp_{\tilde{G}}(B)}^{\tilde{G}}M
    \end{tikzcd}
  \end{equation}
  and
  \begin{equation}\label{poset iso induced by fong 2022-12-12 17:56:15}
    \begin{tikzcd}[row sep=0pt]
      (\stautilt \beta)^{\star \star} \ar[r]&(\stautilt \tilde{B})^{\star \star \star}\\
      M\ar[r,mapsto]&\induc_{\inertiagp_{\tilde{G}}(B)}^{\tilde{G}}M.
    \end{tikzcd}
  \end{equation}
  By \Cref{main-thm-2 block}, we get the poset isomorphism
  \begin{equation}\label{poset iso main 2 block 2022-12-12 17:55:25}
    \begin{tikzcd}[row sep=0pt]
      (\stautilt B)^{\inertiagp_{\tilde{G}}(B)} \ar[r]&(\stautilt \beta)^{\star \star}\\
      M\ar[r,mapsto]&\beta \induc_{G}^{\inertiagp_{\tilde{G}}(B)}M.
    \end{tikzcd}
  \end{equation}
  By \Cref{fong conseq 2022-12-12 17:51:31}, the map \eqref{iso block general 2022-12-15 18:23:37} is the composition of the poset isomorphisms \eqref{poset iso main 2 block 2022-12-12 17:55:25} and \eqref{poset iso induced by fong 2022-12-12 17:56:15}.
  Hence, we complete the proof.
\end{proof}

As an application of \Cref{cor1 block general},
we consider the case that \(\inertiagp_{\tilde{G}}(B)/G\) is a \(p\)-group.
The following \namecref{p-extension block general} is a significant generalization of \cite[Theorem 1.2]{MR4243358} and \cite[Theorem 15]{MR3856858}.

\begin{theorem}\label{p-extension block general}
  Let \(G\) be a normal subgroup of a finite group \(\tilde{G}\), \(B\) a block of \(kG\) and \(\tilde{B}\) a block of \(k\tilde{G}\) covering \(B\).
  If the quotient group \(\inertiagp_{\tilde{G}}(B)/G\) is a
  \(p\)-group, then the functor \(\induc_{G}^{\tilde{G}}\) induces an isomorphism as partially ordered sets between \((\stautilt B)^{\inertiagp_{\tilde{G}}(B)}\) and \(\stautilt \tilde{B}\),
  where \((\stautilt B)^{\inertiagp_{\tilde{G}}(B)}\) is
  the subset of \(\stautilt B\) consisting of \(\inertiagp_{\tilde{G}}(B)\)-invariant support \(\tau\)-tilting \(B\)-modules.
\end{theorem}
\begin{proof}
  It immediately follows from \Cref{covering block p-power 2022-12-12 17:59:37,}, \Cref{main-thm-2 block general} \ref{main-thm-2 block general item3}, \Cref{cor1 block general} and the fact that the only simple \(k(\inertiagp_{\tilde{G}}(B)/G)\)-module is the trivial \(k(\inertiagp_{\tilde{G}}(B)/G)\)-module.
\end{proof}
% \bibliographystyle{abbrvdoi}
% \bibliography{KoshioCite}
% abbrvDOI.bst by David Kotz.
% You must add this to your latex preamble:
%    \usepackage{url}
% or, if you want hyperlinks:
%    \usepackage{hyperref}

\Addresses
\end{document}